\documentclass[12pt]{article}
\usepackage[dvips]{graphicx}
\usepackage{amsmath,amsfonts,amssymb,amsthm,color}

\usepackage{natbib}

\newtheorem{theorem}{Theorem}
\newtheorem{lemma}{Lemma}
\newtheorem{corollary}{Corollary}

\newtheorem{proposition}{Proposition}

\newtheorem{definition}{Definition}
\newtheorem{example}{Example}

\usepackage{bm}
\bmdefine{\Bt}{t}
\bmdefine{\BX}{X}
\bmdefine{\BY}{Y}
\bmdefine{\BZ}{Z}
\bmdefine{\BB}{B}
\bmdefine{\BM}{M}
\bmdefine{\BD}{D}
\bmdefine{\Bi}{i}
\bmdefine{\Bj}{j}
\bmdefine{\Bk}{k}
\bmdefine{\Bx}{x}
\bmdefine{\By}{y}
\bmdefine{\Bz}{z}
\bmdefine{\Bv}{v}
\bmdefine{\Bw}{w}
\bmdefine{\Bn}{n}
\bmdefine{\Ba}{a}
\bmdefine{\Bb}{b}
\bmdefine{\Bc}{c}
\bmdefine{\Be}{e}
\bmdefine{\Bu}{u}
\bmdefine{\Bp}{p}
\bmdefine{\Bzero}{0}
\bmdefine{\Bone}{1}

\newcommand{\R}{{\mathbb R}}
\newcommand{\cB}{{\cal B}}
\newcommand{\cC}{{\cal C}}

\newcommand{\cI}{{\cal I}}
\newcommand{\cJ}{{\cal J}}
\newcommand{\cS}{{\cal S}}

\newcommand{\ci}{\perp\!\!\!\perp}  

\newcommand{\red}{\mathop{{\rm red}}}

\title{Hierarchical subspace models for contingency tables}
\author{
Hisayuki Hara\footnote{
Faculty of Economics, Niigata University},
Tomonari Sei\footnote{
Department of Mathematics,
 Faculty of Science and Technology,
 Keio University} \ \  and 
Akimichi Takemura\footnote{
Graduate School of Information Science and Technology, 
University of Tokyo}
\footnote{CREST, JST}
}
\date{\today}

\begin{document}
\maketitle

\begin{abstract}
 For statistical analysis of multiway contingency tables we propose
 modeling interaction terms in each maximal compact component of a
 hierarchical model.  By this approach we can search for  
 parsimonious models with smaller degrees of freedom than the
 usual hierarchical model, while preserving the localization property of
 the inference in the hierarchical model. 
 This approach also enable us to evaluate the localization property of a
 given log-affine model. 
 We discuss estimation and exacts tests of the proposed model and
 illustrate  the advantage of the proposed modeling with some data 
 sets. 
\end{abstract}


\section{Introduction}
Modeling of the interaction term is an important topic for two-way
contingency tables, because there is a large gap between the
complete independence model and the saturated model.   
This problem is clearly of importance for contingency tables with three 
or more factors.  However modeling strategies of higher order
interaction terms have not been fully discussed in literature.  
In this paper we 
establish a general mathematical framework for
modeling interaction terms of multiway
contingency tables by considering each maximal compact component of a 
hierarchical model.  

For two-way contingency tables the uniform association model
(\citet{goodman1979jasa,goodman1985as}) and the RC association model
(\citet{goodman1979jasa,goodman1985as,kuriki2005}) are often used
for modeling interaction terms.
In the analysis of agreement among raters,  where data are summarized
as square contingency tables with the same categories, many models 
with interaction in diagonal elements and their extension to multiway tables 
have been considered (e.g. \citet{Tanner-Young1985},
\citet{tomizawa09}). 
\citet{hirotsu-1997} proposed a two-way change point model and
\citet{hty-jpaa} generalized it to a subtable sum model. 
For multiway contingency tables \citet{Hojsgaard} considered the split
model as a generalization of graphical models. 
The context specific interaction model defined by \citet{Hojsgaard2004}
is a more general model than the split model.
In this article we give a unified treatment of these models as submodels
of hierarchical models and consider their extension to the models for
higher dimensional tables from viewpoints of decomposition and
conditional independence structure of the models.      

Conditional independence structure of a log-affine model is described by
a graph. 
Such a graph is called an independence graph.
In a usual hierarchical model, the likelihood is factorized to submodels 
induced by each compact components (\cite{Malvestuto-Moscarini}) of the
simplicial complex determining the model. 
By this factorization, statistical inference on a hierarchical model can  
be localized through the decomposition of the simplicial complex for the
model. 
The possibility of localizing the inference of a given hierarchical
model has been well studied by many authors (e.g. \citet{haberman1974},  
\citet{geng1989}, \citet{Malvestuto-Moscarini},
\citet{Badsberg-Malvestuto}, \citet{lauritzen1996}). 

In a usual hierarchical model each maximal interaction effect is
saturated, i.e.\ there is no restriction on the parameters for maximal 
interaction effects.  
However we can consider the modeling for interaction effects of a given
hierarchical model.  
In the modeling process, it is sometimes advantageous to preserve the conditional
independence structure and localization property of the hierarchical
model and to treat each marginal model corresponding to each compact
component of the hierarchical model separately.  
The resulting model is a submodel of the hierarchical model. 
Throughout this paper we assume that the model is log-affine.
When a log-affine model is a submodel of a given hierarchical model, 
the log-affine model has the same conditional independence structure as
the hierarchical model.   
As we will discuss in Section 3, however, the log-affine model does not
necessarily have the same localization property as the hierarchical
model. 
Therefore the localization property of a given log-affine model is not
trivial in general.

In this article we define a hierarchical subspace model by a log-affine
model possessing the same localization property as a given hierarchical
model and discuss the localization property of the log-affine model. 
As pointed out by referees, ideas similar to our hierarchical subspace model
have been discussed in many contexts.
Sociologists have been employing marginal modeling, where a few
important marginals  
are first modeled and they are combined into a joint model.  
\citet{dobra-fienberg-2000pnas} presented maximum likelihood estimation
and bounds for cell entries for reducible models and discuss generalizations
to nongraphical loglinear models.  By our formulation of
the hierarchical subspace model we can discuss these models in a unified
framework.

The organization of the paper is as follows.
In Section \ref{subsec:definition} we give a brief review on
log-affine models and we summarize some basic facts on graphs and hypergraphs.
In Section \ref{sec:model} we define the hierarchical subspace model and
discuss the localization of inference through the decomposition of the
model.  
We show that for a given log-affine model there exists the smallest
decomposable model possessing the same localization property of the
inference. 
In Section \ref{sec:split} we study the split model in the framework of
this paper. 
In Section \ref{sec:markov-basis} we present construction of Markov bases
for conditional tests of our model based on the argument in
\citet{dobra-sullivant} for the hierarchical model.
In Section \ref{sec:examples} we show some real data examples.
Some concluding remarks are given in Section \ref{sec:remarks}.

\section{Definitions and notations}
\label{subsec:definition}
\subsection{Log-affine model and hierarchical model for contingency
  tables}
\label{sec:log-affine}
In this section we summarize basic definitions and notations of
log-affine model and hierarchical model.    
We follow definitions and notations of \citet{darroch-speed} and
\citet{lauritzen1996}.  

Let $V=\R^{I_1\times\dots\times I_m}$ denote the set of 
$I_1\times\dots\times I_m$ tables with real entries, 
where $I_j\geq 2$ for all $j$.
$V$ is considered as an $I_1\times\dots\times I_m$-dimensional real
vector space of functions (tables) from  
${\cal I}=[I_1]\times\dots \times [I_m]$ to $\R$, where $[J]$ denotes
$\{1,\dots,J\}$.   
A probability distribution over $\cal I$ is denoted by 
$\{p(\Bi), \Bi\in {\cal I}\}$.
Let $L$ be a linear subspace of $V$. 
A {\em log-affine model} ${\cal M}(L)$ specified by $L$ is
given by the class of probability functions satisfying 
$\log p(\cdot) \in L$, where $\log p(\cdot)$ denotes the vector   
$\{\log p(\Bi), \Bi\in {\cal I}\}$
(Chapter 4 of \citet{lauritzen1996}).
In the following we only consider linear subspaces of
$V$ containing the constant function 1.

Let $D$ be a subset of $[m]$. $\Bi_D=\{i_j, j\in D\}$  is a $D$-marginal cell. ${\cal
I}_D=\prod_{j\in D}[I_j]$ denotes 
the set of $D$-marginal cells. 
$p(\Bi_D)$ and $x(\Bi_D)$ denote the marginal probability of a
probability distribution $p(\cdot)$ and the marginal frequency of a
contingency table  
$\Bx=\{x(\Bi), \Bi\in {\cal I}\}$, respectively, that is, 
\[
 p(\Bi_D) := \sum_{i_{[m] \setminus D} \in {\cal I}_{[m] \setminus D}}
 p(\Bi),\quad    
 x(\Bi_D) := \sum_{i_{[m] \setminus D} \in {\cal I}_{[m] \setminus D}}
 x(\Bi). 
\]
Define $n := \sum_{i \in {\cal I}} x(\Bi)$, which is the total frequency. 
Denote by $\hat{p}(\Bi)$ and $\hat{p}(\Bi_D)$ the maximum likelihood
estimator (MLE) of $p(\Bi)$ and $p(\Bi_D)$, respectively.
As in \citet{darroch-speed} or \citet{lauritzen1996}, let
\[
F_D = \{ \psi \in V \mid \psi(i_1,\dots,i_m)=\psi(i_1',\dots,i_m')
\ \text{if}\  i_h=i_h', \forall h\in D\}
\]
denote the set of functions depending only on $\Bi_D$.
$F_D$ can be identified with $\R^{I_D}$, where $I_D = \prod_{h \in D} I_h$, 
and especially we note that $F_{[m]} = V$.
For a subspace $L$ of $V$ and $D\subset [m]$, we say that $D$ is {\em
saturated} in $L$ if $F_D \subset L$.  
Then we note the following proposition.
\begin{proposition}
 $D$ is saturated in $L$ if and only if the sufficient statistic
 for ${\cal M}(L)$ fixes all the $D$-marginals of the contingency table.  
\end{proposition}
\begin{proof}
The sufficient statistic for ${\cal M}(L)$ is usually described
by taking a basis of $L$. 
Let $d=\dim L$ and take a basis $\phi_1, \dots, \phi_d$ of $L$.
Then a sufficient statistic for ${\cal M}(L)$ is given as
$\{\sum_{\Bi \in {\cal I}} \phi_j(\Bi)x(\Bi), j=1,\dots,d\}$.
However if we allow redundancy, we can define the
sufficient statistic of $L$ just by 
$\{\sum_{\Bi \in {\cal I}} \phi(\Bi)x(\Bi), \forall \phi(\cdot)\in L\}$.
On the other hand the sufficient statistic for $F_D$ is given
by the set of $D$-marginal frequencies $\{x(\Bi_D), \Bi_D\in {\cal I}_D\}$,
or equivalently by $\{\sum_{\Bi \in {\cal I}} \phi(\Bi)x(\Bi), \forall \phi(\cdot)\in F_D\}$
if we allow redundancy.
Hence the sufficient statistic of $L$ fixes all $x(\Bi_D)$ if and only
if $F_D\subset L$.
\end{proof}
Note that if $D$ is saturated in $L$, then every $E \subset D$ is
saturated in $L$ because $F_E \subset F_D$. 

Let $\Delta$ denote a simplicial complex on $[m]$ and let $\red \Delta$
denote the set of maximal elements, i.e. facets, of $\Delta$ (Chapter 2
of \citet{lauritzen1996}). 
For a subset $D$, define the subcomplex $\Delta(D) := \{D \cap E \mid E \in \Delta\}$. 
The hierarchical model ${\cal M} (H_\Delta)$ associated with $\Delta$ is
defined as 
\[
\log p(\cdot) \in H_\Delta :=
\sum_{D\in \red \Delta} F_D, 
\]
where the right-hand side is the summation
of vector spaces. 
Noting that 
\[
H_\Delta = \left\{
\sum_{D\in \red \Delta} \phi_D(\cdot) \ \mid  \ 
\phi_D(\cdot) \in F_D, D \in \red \Delta
\right\}, 
\]
we have $H_{\Delta \cap \Delta'} = H_\Delta \cap H_{\Delta'}$.

Let $G_\Delta$ be a graph with the vertex set $[m]$ and an
edge between $v, v' \in [m]$ if and only if there exists
$D \in \Delta$ such that $v, v' \in D$. 
Then $G_\Delta$ is called an {\em independence graph} of 
$\Delta$ (\cite{dobra-sullivant}). 
$G_\Delta$ shows an conditional independence structure of 
${\cal M}(H_\Delta)$, i.e., 
if two vertices $v$ and $v'$ are not adjacent each other, the
corresponding variables are conditionally independent given the rest of
variables.
If $\red\Delta$ is the set of maximal cliques of
$G_\Delta$, ${\cal M}(H_\Delta)$ is called a graphical model.
When $G_\Delta$ is chordal, 
a graphical model $H_\Delta$ is called a decomposable model.

\subsection{Basic facts on hypergraphs}
We note that $\red\Delta$ is considered as a hypergraph.
Here we summarize some notions on hypergraphs according to 
\citet{lauritzen1996} and 
\citet{Malvestuto-Moscarini}. 

A hypergraph is reduced if its edges are pairwise inclusion-incomparable
sets. 
Hence $\red\Delta$ is reduced. 
A subset of a hyperedge 
is called a {\em partial edge}. 
A subhypergraph of $\red\Delta$ is a hypergraph whose edges are all
partial edges of $\red\Delta$. 
A subhypergraph of $\red\Delta$ induced by a nonempty subset $E$ of
$[m]$ is $\red\Delta(E)$. 
We note that $\red\Delta(E)$ is a reduced hypergraph whose edges are the
maximal edges of the hypergraph $\{D \cap E \mid D \in \red\Delta\}$. 

Two vertices $v$ and $v'$ are called adjacent in $\red\Delta$ when 
they are also adjacent in $G_\Delta$.
Two vertices $v$ and $v'$ are connected if they are connected in
$G_\Delta$. 
A hypergraph is connected if every pair of two vertices is connected.
A hypergraph is called disconnected if it is not connected.

A partial edge $S$ is a separator of $\red\Delta$ if the subhypergraph
of $\red\Delta$ induced by $[m] \setminus S$ is disconnected. 
For every partial edge separator, there exist three non-empty and
disjoint subsets $\{A,B,S\}$, $A \cup B \cup S =[m]$ satisfying that
$\red\Delta(A)$ and $\red\Delta(B)$ are disconnected.  
Then $\{A,B,S\}$ is called a decomposition of $\red\Delta$. 
For two vertices $u$ and $v$, if there is a decomposition $\{A,B,S\}$
such that $u \in A$  and $v \in B$, we say $S$ separates $u$ and $v$. 
A partial edge separator $S$ of $\red\Delta$ is called a {\em divider}
if there exist two vertices $u,v \in [m]$ that are separated by $S$ but
by no proper subset of $S$. 
If two vertices $u,v \in [m]$ are not separated by any partial edges, 
$u$ and $v$ are called tightly connected.
A subset $C \subset [m]$ is called a {\em compact component} 
if any two vertices in $C$ are tightly connected. 
Denote the set of maximal compact components of $\red\Delta$ by 
$\cal C$. 
Then there exists a sequence of maximal compact components 
$C_1,\ldots,C_{\vert {\cal C} \vert}$ such that 
\[
 (C_1 \cup \cdots \cup C_{k-1}) \cap C_k = S_k
\]
and $S_k$, $k=2,\ldots, \vert {\cal C} \vert$ are dividers of
$\red\Delta$.  
We denote ${\cal S} = \{S_2,\ldots,S_{\vert {\cal C} \vert}\}$.  
${\cal S}$ is a multiset in general.
${\cal C}$ is obtained by decomposing $\red\Delta$ recursively by
dividers. 

By definition it is clear that $v$ and $v'$ are adjacent to each other
in $\red\Delta$ if and only if they are adjacent in $G_\Delta$.
Therefore $\red\Delta$ also gives the conditional independence structure
of the hierarchical model ${\cal M}(H_\Delta)$.
The cell probability $p(\Bi)$ of hierarchical model ${\cal M}(H_\Delta)$ is
factorized as  
\begin{equation}
  p(\bm{i}) = 
  \frac{\prod_{C \in {\cal C}}p(\bm{i}_C)}
  {\prod_{S \in {\cal S}}p(\bm{i}_S)}, 
\end{equation}
where the marginal models $p(\Bi_C)$ and $p(\Bi_S)$ are hierarchical
models ${\cal M}(H_{\Delta(C)})$ and ${\cal M}(H_{\Delta(S)})$, respectively. 
Then the MLE is written as  
\begin{equation}
 \label{eq:cell.prob}
  \hat p(\bm{i}) = 
  \frac{\prod_{C \in {\cal C}}\hat p(\bm{i}_C)}
  {\prod_{S \in {\cal S}}\hat p(\bm{i}_S)} = 
  \frac{\prod_{C \in {\cal C}}\hat p(\bm{i}_C)}
  {\prod_{S \in {\cal S}}x(\bm{i}_S)/n}, 
\end{equation}
and the computation of the MLE is localized to the marginal model
corresponding to each compact component and the localization 
corresponds to the decomposition of $\red\Delta$. 

\begin{example}
 \label{ex:3way}
 \begin{figure}[b]
  \centering
  \includegraphics{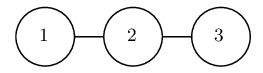}
  \caption{3-way conditional independence model}
  \label{fig:3way}
 \end{figure}
 Consider the decomposable graphical model for three-way contingency
 tables corresponding to the graph in Figure \ref{fig:3way}.  
 The model is described as
 \begin{equation}
  \label{eq:cond-indep}
   \log p(\Bi)=a(i_1, i_2) + b(i_2,i_3).
 \end{equation}
 In this model 
 $\Delta = \{\emptyset,\{1\},\{2\},\{3\},\{1,2\},\{2,3\}\}$ and 
 $\red\Delta = \{\{1,2\},\{2,3\}\}$, respectively, 
 and the corresponding
 linear subspace is $H_\Delta = F_{\{1,2\}} + F_{\{2,3\}}$. 
 We note that $a(i_1,i_2)$'s and $b(i_2,i_3)$'s are
 free parameters. 
 Since the model satisfies $i_1 \ci i_3 \mid i_2$, 
 $p(\Bi)$ is written by 
 \begin{equation}
  \label{eq:p}
   p(\Bi)= \frac{p(\Bi_{\{1,2\}}) p(\Bi_{\{2,3\}})}{p(i_2)}.
 \end{equation}
 The marginal models $p(\Bi_{\{1,2\}})$, $p(\Bi_{\{2,3\}})$ and $p(i_2)$
 are saturated models corresponding to $F_{\{1,2\}}$, 
 $F_{\{2,3\}}$ and $F_{\{2\}}$, respectively. 
 Then the MLE of $p(\Bi)$ is obtained by
 \begin{equation}
  \label{eq:MLE}
   \hat p(\Bi) = \frac{\hat p(\Bi_{\{1,2\}}) \hat p(\Bi_{\{1,2\}})}
   {\hat p (i_2)}
   = \frac{x(\Bi_{\{1,2\}}) x(\Bi_{\{2,3\}})}{n x(i_2)}, 
 \end{equation}
 where $\hat p(\Bi_{\{1,2\}})$, $\hat p(\Bi_{\{2,3\}})$ and $\hat p(i_2)$ are the MLE of 
 $p(\Bi_{\{1,2\}})$, $p(\Bi_{\{2,3\}})$ and $p(i_2)$, respectively. 

 Now consider modeling of two-way interaction terms.
 Suppose that we have 
 known functions $\phi(\Bi_{\{1,2\}})$ depending only on
 $\Bi_{\{1,2\}}=(i_1,i_2)$ and   
 $\psi(\Bi_{\{2,3\}})$ depending only on $\Bi_{\{2,3\}} =(i_2, i_3)$.
 Separating main effects, consider the following submodel of
 (\ref{eq:cond-indep}), 
 \begin{equation}
 \label{eq:3way-ci-sub}
  \log p(\Bi) =\alpha(i_1) + \beta(i_2) + \gamma(i_3) 
  + \delta \phi(\Bi_{\{1,2\}}) + \delta' \psi(\Bi_{\{2,3\}}). 
 \end{equation}
 The model (\ref{eq:cond-indep}) is still log-affine.
 Let $L$ be the linear subspace corresponding to this model. 
 Then $L$ is a linear subspace of $F_\Delta$. 

 The parameters of this model are $\{\alpha(i_1)\}_{i_1=1}^{I_1}, 
 \{\beta(i_2)\}_{i_2=1}^{I_2}, \{\gamma(i_3)\}_{i_3=1}^{I_3}$ and  
 $\delta, \delta'$. 
 The uniform association model 
 is specified  by $\phi(\Bi_{\{1,2\}})= i_1 i_2$.  The 
 change point model
 in \citet{hirotsu-1997} is specified by
 \[
 \phi(\Bi_{\{1,2\}})=
 \begin{cases} 1, & \text{if } i_1 \le I^\prime_1 
  \text{ and } i_2 \le I^\prime_2, \\
  0, & \text{otherwise},
 \end{cases}
 \]
 where $1\le I^\prime_1 < I_1$, $1\le I^\prime_2 < I_2$.  
 Similarly we can  specify $\psi(\Bi_{\{2,3\}})$ according to many well known
 models. 
 
 Since the model (\ref{eq:3way-ci-sub}) is a submodel of the model
 (\ref{eq:cond-indep}),  
 $i_1 \ci i_3 \mid i_2$ still holds for (\ref{eq:3way-ci-sub}) and 
 $p(\Bi)$ is written as (\ref{eq:p}), where 
 we note that the marginal models $p(\Bi_{\{1,2\}})$ and
 $p(\Bi_{\{2,3\}})$ are written by 
 \begin{equation}
  \label{eq:ij-marginal}
   \log p(\Bi_{\{1,2\}}) = \alpha(i_1) + \beta(i_2) 
   + \delta \phi(\Bi_{\{1,2\}})
 \end{equation}
 and 
 \begin{equation}
  \label{eq:jk-marginal}
   \log p(\Bi_{\{2,3\}}) = \beta(i_2) + \gamma(i_3) 
   + \delta^\prime \psi(\Bi_{\{2,3\}}), 
 \end{equation} 
 respectively. 
 Moreover, since $\{\beta(i_2)\}_{i_2=1}^{I_2}$ in 
 (\ref{eq:3way-ci-sub})
 are free parameters, 
 $F_2$ is saturated in $L$.
 Therefore the MLE of $p(\Bi)$ is written by
 \begin{equation}
  \label{eq:ci-MLE}
   \hat p(\Bi) = \frac{\hat p(\Bi_{\{1,2\}}) \hat p(\Bi_{\{1,2\}})}
   {\hat p (i_2)}
   = \frac{\hat p(\Bi_{\{1,2\}}) \hat p(\Bi_{\{2,3\}})}{x(i_2)/n}.
 \end{equation}
 Therefore the maximum likelihood estimation of the model
 (\ref{eq:3way-ci-sub}) is also localized to estimations of two marginal
 models in the same way as the hierarchical model (\ref{eq:MLE}).  

 Note that although we use the same notation for $\beta(i_2)$ in
 (\ref{eq:3way-ci-sub}), (\ref{eq:ij-marginal}) and
 (\ref{eq:jk-marginal}) for simplicity, they are different parameters (as
 functions of cell probabilities). 
 If we distinguish them by $\beta(i_2)^{(123)}$, $\beta(i_2)^{(12)}$,
 $\beta(i_2)^{(23)}$ in (\ref{eq:3way-ci-sub}), (\ref{eq:ij-marginal}),
 (\ref{eq:jk-marginal}), respectively, then they are connected as 
 $\beta(i_2)^{(123)}=\beta(i_2)^{(12)} +  \beta(i_2)^{(23)} -\log p(i_2)$.
 Accordingly, in view of (\ref{eq:ci-MLE}), the maximum likelihood
 estimates are connected as 
 $\hat\beta(i_2)^{(123)}=\hat\beta(i_2)^{(12)} + \hat \beta(i_2)^{(23)} -\log
 (x(i_2)/n)$. 
\end{example}

 When a log-affine model has the same localization property as a given
 hierarchical model as seen in this example, we call the model 
 a hierarchical subspace model of the hierarchical model.   
 Actually the model (\ref{eq:3way-ci-sub}) is a hierarchical subspace
 model of  (\ref{eq:cond-indep}). 
 In the next section we give a precise definition of the hierarchical subspace
 model. 

\section{Hierarchical subspace models and their decompositions}
\label{sec:model}
\subsection{Conformality of log-affine model}
For defining our hierarchical subspace model, 
we introduce the notion of conformality of a hierarchical model.
As an illustrating example, we again consider the three-way conditional
independence model in Example \ref{ex:3way}. 
In (\ref{eq:3way-ci-sub}) it is important to note that $\delta$ and
$\delta'$ are free parameters.   
Now consider the following model imposed an additional constraint $H:
\delta=\delta'$ on (\ref{eq:3way-ci-sub}): 
\begin{equation}
 \label{eq:common-delta}
  \log p(\Bi) = \alpha(i_1) + \beta(i_2) + \gamma(i_3) 
  + \delta (\phi(\Bi_{\{1,2\}}) + \psi(\Bi_{\{2,3\}})). 
\end{equation}
This model is still log-affine and the conditional independence  $i_1 \ci i_3 \mid i_2$ holds. 
However, since $\delta$ is shared by two interaction terms for 
$\Bi_{\{1,2\}}$ and $\Bi_{\{2,3\}}$, both 
$x(\Bi_{\{1,2\}})$ and $x(\Bi_{\{2,3\}})$ are relevant for the
estimation of the common value of $\delta$. 
Therefore we can not localize estimation of the parameters to
two marginal tables.   
We now formulate the above notion of no restriction on parameters across
maximal compact components by defining the notion of conformality of
linear subspaces. 
\begin{definition}
\label{def:conformality}
Let $W_1, \dots, W_K$ be  linear subspaces of $V$. 
A subspace 
$L$ is {\em conformal} to $\{W_j\}_{j=1}^K$ if 
\[
L = (L\cap W_1) + \cdots + (L\cap W_K).
\]
\end{definition}
Any $L$ conformal to $\{W_j\}_{j=1}^K$ is clearly a subspace of 
$W=W_1 + \cdots + W_K$.
Note that if $L$ is a subspace of $W$ then the relation
$L= L\cap W \supset (L\cap W_1) + \cdots + (L\cap W_K)$ always holds but the
inclusion is strict in general.  
We note that $H_\Delta$ satisfies
\begin{equation}
\label{eq:conformality}
H_\Delta = \sum_{C \in {\cal C}} L \cap L_C
\end{equation}
and therefore $H_\Delta$ is conformal to ${\cal C}$.

\begin{example}
\label{ex:conf-3way}
 Consider the models (\ref{eq:3way-ci-sub}) and (\ref{eq:common-delta})
 again. 
 Let $L$ and $L'$ denote the corresponding subspaces of the models
 (\ref{eq:3way-ci-sub}) and (\ref{eq:common-delta}), respectively.
 Let $K = 2$ and let $W_1 := F_{\{1,2\}}$ and $W_2 := F_{\{2,3\}}$.
 In the case of the model (\ref{eq:3way-ci-sub}), 
 \[
 L \cap W_1 = \{\alpha(i_1) + \beta(i_2) + \delta \phi(\Bi_{12})\},
 \quad 
 L \cap W_2 = \{\beta(i_2) + \gamma(i_3) + \delta' \psi(\Bi_{23})\}.
 \]
 Hence $L = (L\cap W_1) + (L\cap W_2)$ 
 is conformal to two marginal spaces $\{F_{\{1,2\}}, F_{\{2,3\}}\}$. 
 In the case of the model (\ref{eq:common-delta}), however, 
 \[
 L' \cap W_1 = \{\alpha(i_1) + \beta(i_2)\}, \quad
 L' \cap W_2 = \{\beta(i_2) + \gamma(i_3)\}.
 \]
 Hence 
 $(L' \cap W_1) + (L' \cap W_2) = \{\alpha(i_1) + \beta(i_2) +
 \gamma(i_3)\}$ and    
 $L'$ is not conformal to $\{F_{\{1,2\}}, F_{\{2,3\}}\}$.
\end{example}
\subsection{Hierarchical subspace model}
We now present the following definition of a hierarchical subspace model.

\begin{definition}
\label{definition:1} 
 Let $\Delta$ be a simplicial complex and $H_\Delta$ be a subspace of
 the corresponding hierarchical model. 
 Then the log-affine model ${\cal M}(L)$ for a subspace $L$ 
 is a hierarchical subspace model (HSM) of $H_\Delta$
 if the following conditions hold:
 \vspace{-5pt}
 \begin{enumerate}
  \setlength{\itemsep}{0pt}
  \item Each divider $S\in {\cal S}$ of $\red\Delta$ is saturated in
	$L$, i.e. 
	$F_S \cap L = F_S$.
  \item $L$ is conformal to the set of subspaces $\{F_C, C\in {\cal C}\}$.
 \end{enumerate}
\end{definition}

By condition 1 of HSM the conditional independence structure of
$H_\Delta$ is preserved in $L$.  
Condition 2 together with condition 1 guarantees that the statistical inference is
localized to each $C$.

On the computation of the MLE we can generalize  (\ref{eq:cell.prob})
to HSM as follows.

\begin{theorem}
 \label{thm:mle-local}
 The MLE $\hat{p}(\bm{i})$ of cell probabilities for HSM 
 of $H_\Delta$ satisfies 
 \begin{equation}
  \label{eq:mle}
   \hat{p}(\bm{i}) = 
   \frac{\prod_{C \in {\cal C}}\hat{p}(\bm{i}_C)}
   {\prod_{S \in {\cal S}}\hat{p}(\bm{i}_S)}
   =
   \frac{\prod_{C \in {\cal C}}\hat{p}(\bm{i}_C)}
   {\prod_{S \in {\cal S}}x(\bm{i}_S)/n}, 
 \end{equation}
 where $\hat{p}(\bm{i}_C)$ coincides with the MLE of the model
 associated with the linear space $L \cap F_C$, which is computed only
 on the marginal table  $x(\bm{i}_C)$.
\end{theorem}

\begin{proof}
 By induction on the number of compact components $|{\cal C}|$ of 
 $\red\Delta$, it is
 sufficient to consider the case $\cC=\{C_1,C_2\}$ with $S=C_1 \cap C_2$. 
 The MLE of the model ${\cal M}(L)$ is the maximizer of 
 $\sum_{\Bi}x(\Bi) \log p(\Bi)$ subject to $\log p(\cdot)\in L$ and 
 $\sum_{\Bi}p(\Bi)=1$.
 By Condition 2 we write $\log p(\cdot)=\theta_{C_1}+\theta_{C_2}$ with  
 $\theta_{C_1} \in L \cap F_{C_1}$ and $\theta_{C_2} \in L \cap F_{C_2}$.
 Since $F_S$ is saturated both in $L\cap F_{C_1}$ and 
 $L\cap F_{C_2}$, we can assume  
 $\sum_{\Bi_{C_1 \setminus S}}e^{\theta_{C_1}(\Bi_{C_1})}=1$ for each 
 $\Bi_S$ without loss of generality. 
 Hence the problem is decomposed into two parts: maximization of
 $\sum_{\Bi_{C_1}} x(\Bi_{C_1})\theta_{C_1}(\Bi_{C_1})$ subject to 
 $\theta_{C_1} \in L\cap F_{C_1}$ and 
 $\sum_{\Bi_{C_1 \setminus S}}e^{\theta_{C_1}(\Bi_{C_1})}=1$, 
 and maximization of 
 $\sum_{\Bi_{C_2}}x(\Bi_{C_2})\theta_{C_2}(\Bi_{C_2})$ subject to
 $\theta_{C_2}\in L\cap F_{C_2}$ and 
 $\sum_{\Bi_{C_2}}e^{\theta_{C_2}(\Bi_{C_2})}=1$.
 Since the maximizer $\hat{\theta}_{C_1}$ does not depend on $C_2$,
 it is computed from the case $C_2=S$. 
 We have 
 $\hat{\theta}_{C_1}(\Bi_{C_1})=\log\{\hat{p}(\Bi_{C_1})/(x(\Bi_S)/n)\}$,
 where $\hat{p}(\Bi_{C_1})$ is the MLE of the model ${\cal M}(L\cap F_{C_1})$.
\end{proof}

This Theorem shows that the computation of the MLE of an HSM of
$H_{\Delta}$ is localized to each $C \in {\cal C}$.  
We note that Theorem \ref{thm:mle-local} depends on Condition 1. 
Even if Condition 1 is not satisfied, the conditional independence
structure of ${\cal M}(H_\Delta)$ is preserved. 
But $\hat{p}(\bm{i}_C)$ is not necessarily the MLE for the marginal model 
${\cal M}(L \cap F_C)$. 

\begin{example}
 \label{ex:HSM}
 By following the argument in Example \ref{ex:conf-3way}, we can easily
 show that the model (\ref{eq:3way-ci-sub}) is an HSM of
 (\ref{eq:cond-indep}).
 On the other hand, since the model (\ref{eq:common-delta}) is not
 conformal to $F_{\{1,2\}}$ and $F_{\{2,3\}}$, 
 the model (\ref{eq:common-delta}) is not an HSM of
 (\ref{eq:cond-indep}). 
 Although the model (\ref{eq:common-delta}) has the same conditional
 independence structure $i_1 \ci i_3 \mid i_2$ depicted in the
 graph in Figure \ref{fig:3way}, 
 the inference is not localized in the same way as the decomposition of
 the graph. 

 As seen in this example, we note that 
 even if a given log-affine model ${\cal M}(L)$ is a subset of a hierarchical 
 model ${\cal M}(H_\Delta)$, the localization property of ${\cal M}(H_\Delta)$ is not  
 necessarily preserved in $L$.  

 However we note that the model (\ref{eq:common-delta}) is an HSM of the
 three-way saturated model. 
 In the saturated model, $\red\Delta = {\cal C} = [m]$ and 
 there is no divider in $\red\Delta$.
 Therefore every log-affine model is an HSM of the saturated model.
 This also means that every log-affine model ${\cal M}(L)$ has a hierarchical
 model for which ${\cal M}(L)$ is an HSM. 
\end{example}

\subsection{Ambient decomposable model of a log-affine model}
\label{subsec:ambient-decomposable}
Suppose that a conditional independence structure of the model is given
by a hypergraph $\red\Delta$. 
By following Definition \ref{definition:1}, 
we can formulate an HSM of $H_\Delta$ by modeling interaction terms 
$L \cap F_D$, $D \in \red \Delta$, under the conditions 
of conformality (\ref{eq:conformality}) and $F_S \subset L$, $S \in \cS$.
Then the resulting model preserves the same localization property as $H_\Delta$.  

Since every log-affine model ${\cal M}(L)$ has a hierarchical model for which
${\cal M}(L)$ is an HSM, a next natural question is to look for a small simplicial  
complex $\Delta$ such that ${\cal M}(L)$ is an HSM of $H_\Delta$.
As mentioned in Example \ref{ex:HSM}, 
even if $L \subset H_\Delta$, the localization property of 
${\cal M}(L)$ does not necessarily correspond to the decomposition of
$\red\Delta$.  
Therefore the question is not trivial. 
We will show in Theorem \ref{thm:acycle} below that for each
log-affine model ${\cal M}(L)$ there exists a natural smallest
decomposable model ${\cal M}(H_{\cal H})$ with respect to inclusion relation, 
such that ${\cal M}(L)$ is an HSM of $H_{\cal H}$.   
Here ${\cal H}$ is the hypergraph corresponding to the decomposable model.
We call such ${\cal M}(H_{\cal H})$ the {\em ambient decomposable model} of 
${\cal M}(L)$.   
The notion of ambient decomposable model is also interpreted as a
classification of log-affine models in terms of decomposition of the
models. 

In order to define the ambient decomposable model, we first introduce
the notion of connectedness and decomposition of a subspace
$L$ separately from those of hypergraphs. 
$L$ is called disconnected if there exists a non-empty proper subset $A$ 
of $[m]$ such that $L$ is conformal to $\{F_A, F_{A^C}\}$, where $A^C$ 
denotes the complement of $A$ in $[m]$. 
We call $L$ {\em connected} if $L$ is not disconnected.  
Now we note the following proposition. 
\begin{proposition}
 When $L$ is disconnected, the variables in $A$ and the variables in
 $A^C$ are independent. 
\end{proposition}
\begin{proof}
 $L= (L \cap F_A) + (L \cap F_{A^C})$
 means that ${\cal M}(L)$ is described as
 $\log p(\Bi) = \phi(\Bi_A) + \psi(\Bi_{A^C})$, where
 $\phi(\cdot) \in F_A$ and $\psi(\cdot) \in F_{A^C}$.
 Therefore $A$ and $A^C$ are independent. 
\end{proof}
Under this definition $L$ can be decomposed into its connected
components. 
By the above proposition, variables in different connected components are
independent. 
Therefore they can be independently modeled in $L$ and can be investigated  
separately.   
Therefore from now on we assume that $L$ is connected.

We need to generalize the notion of partial edge separator of a
hypergraph to our setting. 

\begin{definition}
\label{definition:divider-for-L}
 For a subspace $L$, a non-empty subset $S$ of $[m]$ is called an 
 $L$-separator 
 if $[m]$ is partitioned into three non-empty and
 disjoint subsets $\{A_1,A_2,S\}$ such that 
 \vspace{-5pt}
 \begin{enumerate}
  \setlength{\itemsep}{0pt}
  \item $S$ is saturated in $L$.
  \item  $L$ is conformal to $\{F_{A_1 \cup S}, F_{A_2 \cup S}\}$.
 \end{enumerate}
 Then we call the triple $(A_1,A_2,S)$ a {\em decomposition} of $L$. 
 When the subspace $L$ has a 
 $L$-separator, we call $L$ {\em reducible}. 
 A pair of vertices $v$ and $v'$ are called {\em tightly connected} in
 $L$ if there does not exist a decomposition $(A_1,A_2,S)$ of $L$ such
 that $v \in A_1$ and $v' \in A_2$. 
 When $L$ is not reducible, we call $L$ {\em prime}.
\end{definition}



A set of vertices such that any two of them are tightly
connected in $L$ is called an {\em extended compact component} of $L$.  
We note that 
the notions of $L$-separator, tight connectivity in $L$ and extended compact 
component for a hierarchical model ${\cal M}(H_\Delta)$ are exactly 
the same as the notions of partial edge separator, tight connectivity 
and compact component of the hypergraph $\red\Delta$.  

The set of maximal extended compact components of $L$ is also considered
as a hypergraph and we denote it by ${\cal H}$.  
Denote by $H_{\cal H}$ the subspace of the 
hierarchical model induced by ${\cal H}$. 
Then we have the following theorem.   

\begin{theorem}
 \label{thm:acycle}
 ${\cal M}(H_{\cal H})$ is the smallest decomposable model with respect to
 inclusion relation such that ${\cal M}(L)$ is 
 an HSM of $H_{\cal H}$. 
\end{theorem}

The  following corollary is obvious from (\ref{eq:mle}).  
\begin{corollary}
\label{cor:MLE} 
 The MLE $\hat p(\Bi)$ satisfies
\[
\hat p(\Bi)=\frac{\prod_{C\in {\cal H}} \hat p(\Bi_C)}{
\prod_{S\in {\cal S}} x(\bm{i}_S)/n}, 
\]
where ${\cal S}$ is the set of dividers of ${\cal H}$ 
and $\hat{p}(\Bi_C)$ depends only on the marginal table $x(\Bi_C)$.
\end{corollary}

\bigskip

The rest of this subsection is devoted to a proof of Theorem \ref{thm:acycle}. 
Before we give the proof, we present some lemmas required to prove the
theorem. 

\begin{lemma}
 \label{lemma:partial_edge_separator}
 If $S$ is a 
 $L$-separator, 
 $S$ is also a partial edge separator of the hypergraph ${\cal H}$.
\end{lemma}

\begin{proof}
 Since $S$ is saturated in $L$, 
 $S$ is an extended compact component.
 Hence $S$ is a partial edge of ${\cal H}$. 
 Denote by ${\cal H}([m] \setminus S)$ the subhypergraph of
 ${\cal H}$ induced by $[m] \setminus S$.
 Assume that $S$ is not a separator of ${\cal H}$.
 Then ${\cal H}([m] \setminus S)$ is connected. 
 
 Since $S$ is a separator of $L$, there exists a decomposition
 $(A, B, S)$ of $L$ by definition.
 Define $\tilde{\cal H}(A)$ and $\tilde{\cal H}(B)$ by 
 \[
 \tilde{\cal H}(A) := \{C \in {\cal H} \mid A \cap C \ne \emptyset\}, 
 \quad
 \tilde{\cal H}(B) := \{C \in {\cal H} \mid B \cap C \ne \emptyset\}.
 \]
 Then we have 
 $\tilde{\cal H}(A) \cap \tilde{\cal H}(B) = \emptyset$ which
 contradicts the fact that ${\cal H}([m] \setminus S)$ is connected. 
\end{proof}

When there exists a chordal graph whose set of maximal clique is 
${\cal H}$, ${\cal H}$ is called acyclic. 
By using Lemma \ref{lemma:partial_edge_separator}, 
we can prove the following lemma in the same way as Theorem 5 in 
\citet{Malvestuto-Moscarini}.

\begin{lemma}
 \label{lemma:acyclic}
 ${\cal H}$ is acyclic.
\end{lemma}

Denote by ${\cal S}$ the set of dividers of ${\cal H}$. 

\begin{lemma}
 \label{lemma:divider}
 Suppose $S \in {\cal S}$ is a divider of ${\cal H}$ with a
 decomposition $(A,B,S)$.
 Then $S$ is an $L$-separator 
 with a decomposition
 $(A,B,S)$. 
\end{lemma}

\begin{proof}
 Since $S$ is a divider, there exists a pair of vertices $\{u,v\}$
 such that $S$ is the unique minimal partial edge separating $u$ and $v$. 
 Then there exists a decomposition  $(A,B,S)$ such that 
 $u \in A$ and $v \in B$.   
 Any vertices in $A$ and any vertices in $B$ are not tightly connected
 in $L$.
 This implies that there exists 
 an $L$-separator $S' \subset S$ 
 and a decomposition $(A',B',S')$ of $L$ satisfying 
 $A' \supset A$ and $B' \supset B$.
 From Lemma \ref{lemma:partial_edge_separator}, 
 $S'$ is also a partial edge separator of ${\cal H}$. 
 Noting that $S$ is the unique minimal partial edge of ${\cal H}$
 separating $u$ and $v$, we have $S'=S$.
 Then  $(A,B,S)$ is a decomposition of $L$. 
\end{proof}

Now we provide a proof of Theorem \ref{thm:acycle}. 

\begin{proof}[Proof of Theorem \ref{thm:acycle}]
 It is obvious that $L \subset H_{\cal H}$. 
 From Lemma \ref{lemma:divider}, every divider $S \in {\cal S}$ of 
 ${\cal H}$ is 
 an $L$-separator and hence saturated in $L$.  
 From Lemma \ref{lemma:acyclic}, ${\cal H}$ is considered as the set of
 maximal cliques of a chordal graph ${\cal G}^{\cal H}$. 
 Let $C_k$, $k=1,\ldots,K$, be a perfect sequence of maximal cliques in
 ${\cal G}^{\cal H}$ 
 (see e.g. Section 2.1.3 of \citet{lauritzen1996}). 
 Let 
 \[
 B_k := C_1 \cup C_2 \cup \cdots \cup C_k, \quad 
 R_k := (C_K \cup C_{K-1} \cup \cdots \cup C_k) \setminus S_k, \quad
 S_k := B_{k-1} \cap C_k. 
 \]
 It is known that $S_K$ is a divider of ${\cal H}$ with a decomposition 
 $(B_{K-1}, R_K, S_K)$. 
 From Lemma \ref{lemma:divider}, $S_K$ is 
 an $L$-separator with the same decomposition. 
 Hence $L$ is conformal to $\{F_{B_{K-1}}, F_{C_K}\}$, i.e.
 \[
  L = (L \cap F_{B_{K-1}}) + (L \cap F_{C_K}).
 \]
 In the same way $S_{K-1}$ is 
 an $L$-separator with a decomposition  $(B_{K-2}, R_{K-1}, S_{K-1})$
 and hence $L$ is conformal to $\{F_{B_{K-2}}, F_{C_K \cup C_{K-1}}\}$,
 i.e. 
 \begin{align*}
  L & = (L \cap F_{B_{K-2}}) + (L \cap F_{C_K \cup C_{K-1}})\\
    & = 
  \left[
  \left(
  (L \cap F_{B_{K-1}}) + (L \cap F_{C_K})
  \right) \cap F_{B_{K-2}}
  \right]\\
  & \qquad + 
  \left[
  \left(
  (L \cap F_{B_{K-1}}) + 
  (L \cap F_{C_K})
  \right) \cap F_{C_{K-1} \cup C_K}
  \right]\\
  & = (L \cap F_{B_{K-2}}) + (L \cap F_{C_{K-1}}) +
  (L \cap F_{C_{K}}).
 \end{align*}
 By iterating this procedure, we can obtain
 $L = (L \cap F_{C_1}) + \cdots + (L \cap F_{C_K})$.
 Hence $L$ is conformal to $\{F_C, C \in {\cal H}\}$. 
 Therefore ${\cal M}(L)$ is 
 an HSM of $H_{\cal H}$. 
 
 Suppose that there exists a smaller decomposable model 
 associated with a subspace $F_{\cal H'} \subset H_{\cal H}$ 
 for which ${\cal M}(L)$ is 
 an HSM. 
 Then there exist $C \in {\cal H}$ and a divider $S'$ of ${\cal H}'$ such 
 that $S' \subset C$. 
 This contradicts the fact that any vertices in $C$ are tightly
 connected in $L$. 
\end{proof}




\subsection{Hierarchical models containing a log-affine model}
\label{subsec:closure}

In Theorem \ref{thm:acycle} we have shown the existence
of the smallest decomposable model containing a log-affine model.
Then a natural question is to ask whether there exists 
a smallest hierarchical model with respect to inclusion relation
containing a log-affine model as an HSM. 
In general this does not hold and  we  here discuss  properties of
hierarchical models containing a log-affine model.

As an example consider the model (\ref{eq:common-delta}) again.  
As seen in Example \ref{ex:HSM}, 
(\ref{eq:common-delta}) is a submodel of (\ref{eq:cond-indep})
but is not an HSM of (\ref{eq:cond-indep}). 
The difficulty lies in the fact that a hierarchical model containing $L$
may have a partial edge separator which is not 
an $L$-separator. 

%

Given a subspace $L$
consider the subspace of hierarchical models $H_\Delta$ containing 
$L$: $\{ H_\Delta \mid H_\Delta \supset L \}$.
As mentioned in Section \ref{sec:log-affine}, 
$H_\Delta \cap H_{\Delta'}=H_{\Delta \cap \Delta'}$. 
It follows that there exists the smallest hierarchical model in 
$\{ {\cal M}(H_\Delta) \mid H_\Delta \supset L \}$.  
We call the smallest hierarchical model containing $L$ as 
{\em hierarchical closure} of $L$ and
denote the corresponding simplicial complex and the 
subspace by $\bar \Delta(L)$ and $H_{\bar \Delta(L)}$, respectively. 
Note that for both (\ref{eq:3way-ci-sub}) and (\ref{eq:common-delta}), 
the hierarchical closure is the three-way conditional independence model
(\ref{eq:cond-indep}).
We note that $L$ does not necessarily satisfy the conformality with
respect to the linear subspaces for $\red\bar\Delta(L)$.  
We call ${\cal M}(L)$ a tight hierarchical subspace model if ${\cal M}(L)$ is 
an HSM of $H_{\bar \Delta(L)}$.
If ${\cal M}(L)$ is a tight 
HSM, obviously $\bar \Delta(L)$ is the smallest simplicial complex such
that ${\cal M}(L)$ is its HSM of $H_{\bar \Delta(L)}$.

\begin{figure}[th]
\begin{center}
  \includegraphics[width=0.7\textwidth]{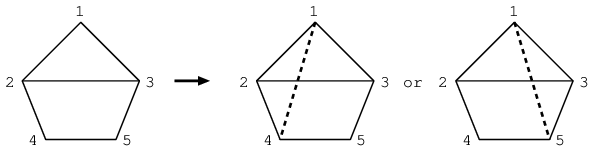}
 \caption{Two ways to cross a divider of the hierarchical closure}
 \label{fig:5wax-ex}
\end{center}
\end{figure}

We now present an example of a log-affine model $L$ of a 5-way
contingency table, which has two minimal hierarchical models
${\cal M}(H_{\Delta_1})$, ${\cal M}(H_{\Delta_2})$, such that ${\cal M}(L)$ is 
an HSM of both of them. 
Consider the following model ${\cal M}(L)$ of 5-way contingency tables:
\begin{align*}
\log p(i_1, \dots, i_5)&=\sum_{j=1}^5 \alpha_{\{j\}}(i_j) + \theta \big(
\psi_{\{1,2\}}(i_1, i_2)+
\psi_{\{1,3\}}(i_1, i_3)+
\psi_{\{2,3\}}(i_2, i_3) \\
& \qquad\qquad +\psi_{\{2,4\}}(i_2, i_4)+
\psi_{\{3,5\}}(i_3, i_5)+
\psi_{\{4,5\}}(i_4, i_5)\big),
\end{align*}
where the main effects  $\alpha_{\{j\}}$'s and $\theta$ are parameters 
and $\psi_{\{j,j'\}}$'s  are fixed functions.  The set of facets of
$\bar \Delta(L)$ is given by
\[
\red \bar \Delta(L) = 
\{\{1,2\}, \{1,3\}, \{2,3\}, \{2,4\}, \{3,5\}, \{4,5\}\}, 
\]
which has a divider $\{2,3\}$.  
On the other hand, since $\psi_{\{2,3\}}(\cdot)$ is a fixed function, 
$L \cap F_{\{2,3\}}$ is not saturated in $L$ and hence 
$\{2,3\}$ is not an $L$-separator. 
Therefore ${\cal M}(L)$ is not an HSM of $H_{\bar\Delta(L)}$ and is not
tight.  
Note that ${\cal M}(L)$ is an HSM of any $H_\Delta$, such that 
$H_\Delta$ does not possess 
a partial edge separator and $L \subset H_\Delta$.
As in Figure \ref{fig:5wax-ex} define
\[
\red \Delta_1= \red\bar\Delta(L) \cup \{\{1,4\}\}, \qquad
\red \Delta_2= \red\bar\Delta(L) \cup \{\{1,5\}\}.
\]
Then ${\cal M}(L)$ is 
an HSM of both $H_{\Delta_1}$ and $H_{\Delta_2}$.

\section{Split model as a hierarchical subspace model}
\label{sec:split}

In this section we give a brief review on the split model by 
\citet{Hojsgaard}.
We first define the context specific interaction (CSI) model
(\citet{Hojsgaard2004}). 
The split model is a particular case of the CSI model.
Recall that $V=\R^{|\cI|}$ is the set of all tables.
For any subset $B$ of $[m]$ and $\Bj_B\in\cI_B$, we
consider a subspace $F^{\Bj_B}$ of $V$ in which only the $\Bj_B$-slice
has nonzero components, that is,
\begin{align*}
 F^{\Bj_B}
 \ &=\ 
 \left\{\psi\in V\mid
 \psi(\Bi)=0\ {\rm if}\ \Bi_B\neq \Bj_B
 \right\}.
 \\
 \ &=\ 
 \left\{\psi\in V\mid\psi(\Bi)
 =f(\Bi_{[m]\setminus B})1_{\{\Bi_B=\Bj_B\}},\ f:\cI_{[m]\setminus B}\to\R\right\}.
\end{align*}
If $B$ is empty, we define $F^{\Bj_{\emptyset}}=V$ with a dummy symbol $\Bj_{\emptyset}$.
For any subsets $B$ and $D$ of $[m]$ and any level $\Bj_B\in \cI_B$,
we define a subspace
\begin{align*}
 F_D^{\Bj_B}\ =\ F_{D\cup B}\cap F^{\Bj_B}
 \ =\ \left\{
 \psi\in V\mid
 \psi(\Bi)=f(\Bi_{D\setminus B})1_{\{\Bi_B=\Bj_B\}},\ 
 f:\cI_{D\setminus B}\to \R
 \right\}.
\end{align*}
The subspace $F_D^{\Bj_B}$ represents {\em a context specific interaction},
that is, an interaction over $\Bi_D$ exists
only if $\Bi_B=\Bj_B$.
The following relation is easily proved:
\begin{align}
 \label{eqn:CSI-easy}
 F_{D\cup B}\ =\ \sum_{\Bj_B\in\cI_B}F_D^{\Bj_B}.
\end{align}
{\em A context specific interaction (CSI) model}
is a direct sum of subspaces $F_D^{\Bj_B}$
for a set of $(\Bj_B,D)$'s.
It is easily shown that any hierarchical model is a CSI model.

Next we define split models.
In order to clarify the definition,
we consider a more general model,
the split subspace model.
The split model is a particular case of the split subspace models.
Although \citet{Hojsgaard} defined the split model on the basis of a graphical model,
we let the graphical model be a decomposable model for simplicity.

Consider a decomposable model ${\cal M}(H_{\Delta})$
with the set of maximal cliques $\cC$.
For each $C\in\cC$ choose a subset $Z(C)\subset C$.
We admit the case where $Z(C)$ is empty.
For each $\Bj_{Z(C)}\in\cI_{Z(C)}$,
choose a subspace $N_C^{\Bj_{Z(C)}}\subset F_C^{\Bj_{Z(C)}}$
such that
\begin{align}
 \label{eqn:CSI-saturated}
 &\forall C'\in \cC\setminus\{C\},\quad
 F_{C\cap C'}^{\Bj_{Z(C)}}\subset N_{C}^{\Bj_{Z(C)}}\subset F_C^{\Bj_{Z(C)}}.
\end{align}
Then a log-affine model ${\cal M}(L)$ is defined by
\begin{align}
 \label{eqn:split-hierarchical}
 L\ =\ \sum_{C\in\cC}N_C,
 \quad N_C\ =\ \sum_{\Bj_{Z(C)}\in\cI_{Z(C)}}
 N_C^{\Bj_{Z(C)}}.
\end{align}
We call ${\cal M}(L)$ {\em a split subspace model with root $\cC$}
if $L$ satisfies (\ref{eqn:CSI-saturated}) and (\ref{eqn:split-hierarchical}).
The following proposition holds.

\begin{proposition}
 \label{lem:split-hierarchical}
 Let ${\cal M}(H_{\Delta})$ be a decomposable model
 with the cliques $\cC$.
 Then any split subspace model ${\cal M}(L)$ with root $\cC$
 is 
 an HSM of $H_{\Delta}$.
\end{proposition}

\begin{proof}
First we prove that $F_S\subset L$ for any divider $S$.
From the definition of dividers of decomposable models,
there exist two cliques $C$ and $C'$ ($C\neq C'$)
such that $S=C'\cap C$.
By the relations (\ref{eqn:CSI-easy}) and (\ref{eqn:CSI-saturated}),
we have 
\begin{align*}
 F_S
 \ \subset\ F_{(C'\cap C)\cup Z(C)}
 \ =\ \sum_{\Bj_{Z(C)}\in\cI_{Z(C)}}F_{C'\cap C}^{\Bj_{Z(C)}}
 \ \subset\ \sum_{\Bj_{Z(C)}\in\cI_{Z(C)}}N_C^{\Bj_{Z(C)}}
 \ =\ N_C.
\end{align*}
Therefore $F_S\subset L$.
Next, we prove that $L$ is conformal to $\{F_C\mid C\in\cC\}$.
Note that 
$N_C^{\Bj_{Z(C)}}\subset F_C^{\Bj_{Z(C)}}\subset F_C$
for any $\Bj_{Z(C)}$ and we have 
$N_C\subset F_{C}$ for each $C\in\cC$.
Since $N_C$ is also a subspace of $L$,
we obtain $N_C\subset L\cap F_C$
and therefore $L=\sum_{C\in\cC}N_C\subset\sum_{C\in\cC}(L\cap F_C)$.
The opposite inclusion is obvious.
\end{proof}


Now we define a split model as a special case of split subspace models.
We say that any decomposable model
is a split model of degree zero.
Then {\em a split model of degree one} is defined as
the decomposition (\ref{eqn:split-hierarchical}) with
\begin{align*}
 N_C^{\Bj_{Z(C)}}
 \ =\ \sum_{D\in \cC_C^{\Bj_{Z(C)}}} F_D^{\Bj_{Z(C)}},
\end{align*}
where $\cC_C^{\Bj_{Z(C)}}$ is
a decomposable model with the vertex set $C\setminus Z(C)$.
Here we assume
\begin{align}
 \label{eqn:split-saturated}
 \forall C'\in\cC\setminus\{C\},\ \exists D\in\cC_C^{\Bj_{Z(C)}}
 \ {\rm s.t.}\ (C\cap C')\setminus Z(C)\subset D
\end{align}
to assure the condition (\ref{eqn:CSI-saturated}).
Split models of degree greater than one
are defined recursively.
See \citet{Hojsgaard} for details.

In Section~\ref{sec:examples}, we will consider
an example of the split model (of degree one).
The following elementary lemma is useful
to obtain the MLE of split models.

\begin{lemma}
 \label{lemma:split-Markov}
 Let $\cI=\bigcup_{\lambda}\cJ_{\lambda}$ be a partition of $\cI$
 and consider subspaces $N_{\lambda}\subset V$ such that
 \begin{align*}
   N_{\lambda}\subset \{\psi\in V\mid \psi(\Bi)=0\ {\rm if}\ \Bi\notin \cJ_{\lambda}\}.
 \end{align*} 
 Then the MLE of the model associated with the subspace $\sum_{\lambda}N_{\lambda}$ 
 is given by  
 $\hat{p}(\Bi)=\sum_{\lambda}(n_{\lambda}/n)\hat{p}_{\lambda}(\Bi)
 1_{\{\Bi\in \cJ_{\lambda}\}}$, where
 $\hat{p}_{\lambda}(\Bi)$ is the MLE of the model ${\cal M}(N_{\lambda})$
 with the total frequency $n_{\lambda}=\sum_{\Bi\in\cI_{\lambda}}x(\Bi)$.
\end{lemma}

\section{Conditional tests of hierarchical subspace models via Markov bases}
\label{sec:markov-basis}
So far we have discussed the localization of the computation of the MLE
for the log-affine model. 
In the hierarchical model, \citet{dobra-sullivant} showed that the
computation of Markov bases is also localized to the computation of the
Markov bases of the marginal model corresponding to each maximal compact
component.  
In this section we generalize the argument to 
an HSM. 

In this section we first give a brief review on Markov bases and
conditional tests based on Markov basis methodology 
(\citet{diaconis-sturmfels}).
Next we generalize the argument of \citet{dobra-sullivant} to the HSM. 

\subsection{Markov basis and conditional test}

Let $\Bb$ be the set of sufficient statistics for ${\cal M}(L)$. 
We assume that the elements of $\Bb$ are integer combinations of the
frequencies $x(\Bi)$.
For a hierarchical model ${\cal M}(H_{\Delta})$, $\Bb$ is written by 
\[
\Bb = \{ x(\bm{i}_{D}), \bm{i}_{D}\in {\cal  I}_D, D \in  \red{\Delta}\}.
\]
We consider $\Bb$ as a column vector with dimension $\nu$.
 
We order the elements of a contingency table $\Bx$ 
lexicographically and consider $\Bx$ as a column vector.  
Then the relation between the joint frequencies $\Bx$ and the marginal
frequencies $\Bb$ is written simply as
\[
\Bb = A \Bx ,
\]
where $A$ is a $\nu \times |{\cal I}|$ integer matrix. 
$A$ is called the configuration for ${\cal M}(L)$.

The conditional distribution of $\Bx$ given  $\Bb$ is exactly a
hypergeometric distribution. 
Usually the goodness of fit of the model is assessed by large sample
approximation. However when the sample size is not large, it is
desirable to use conditional tests based on the exact distribution of
test statistics.
Given $\Bb$, the set
\[
{\cal F}_{\Bb} = \{ \Bx  \ge 0   \mid  \Bb = A \Bx \}
\]
of contingency tables sharing the same $\Bb$ is called a {\em fiber}. 
If we can enumerate all the elements of the fiber which $\Bx$ belongs
to, we can evaluate the null distribution of a test statistic exactly
based on the conditional hypergeometric distribution of $\Bx$. 
However since the number of elements of fibers is too large in general, 
it is difficult to evaluate the null distribution of a test statistic by
the enumeration of elements of a fiber. 

An integer array $\Bz = \{z(\bm{i})\}_{\bm{i} \in {\cal I}}$ of the same
dimension as $\Bx$ is called a {\em move} if $A\Bz = \Bzero$. 
A move is expressed as a difference of its positive part and 
negative part $\Bz = \Bz^+ - \Bz^-$, where $\Bz^+$ and $\Bz^-$ are two
contingency tables in the same fiber. 
We denote a move $\Bz$
\begin{equation}
 \label{eq:moves}
 \Bz = 
[\{\bm{i}_1,\dots,\bm{i}_d\} \Vert 
\{\bm{i}_1^\prime,\dots,\bm{i}_d^\prime\}],
\end{equation}
where $\bm{i}_1,\dots,\bm{i}_d \in {\cal I}$ are cells (with
replication) of positive elements of $\Bz^+$ and
$\bm{i}_1^\prime,\dots,\bm{i}_d^\prime \in {\cal I}$ are cells 
of positive elements of $\Bz^-$. 
$d$ is the sample size of $\Bz^+$ (or $\Bz^-$) and is called a degree of
$\Bz$. 

\begin{example}
 \label{cdem}
 Consider a $3 \times 3$ common diagonal effect model discussed in 
 \citet{hara-takemura-yoshida-diagonal}, 
 \begin{equation}
  \label{eq:cdem}
 \log p(\Bi) = \alpha(i_1) + \beta(i_2) + \delta \phi(\Bi), 
 \end{equation}
 where 
 \begin{equation}
  \label{eq:phi-cdem}
 \phi(\Bi) = \left\{
 \begin{array}{ll}
  1 & i_1 = i_2, \\
  0 & \text{otherwise}.
 \end{array}
 \right.
 \end{equation}
 The sufficient statistic $\Bb$ of this model is the set of row sums,
 column sums and diagonal sums, 
 \[
 \Bb = \left\{
 x(i_1), i_1 \in \{1,2,3\}, \;  x(i_2), i_2 \in \{1,2,3\}, \;
 \sum_{i_1=1}^3 x(i_1 i_1)  \right\}.
 \]
 Then an integer array
 \begin{equation}
  \label{eq:move-array}
  \Bz := 
 \begin{array}{r|r|r|r|}
    \multicolumn{1}{c}{~} &
  \multicolumn{3}{c}{i_2} \\ \cline{2-4}
  & 0 & 1 & -1\\ \cline{2-4}
  i_1 & -1 & 0 & 1\\ \cline{2-4}
  & 1 & -1 & 0\\ \cline{2-4}
 \end{array}
 =
 \begin{array}{|r|r|r|}
  \multicolumn{3}{c}{} \\ \hline
  0 & 1 & 0\\ \hline
  0 & 0 & 1\\ \hline
  1 & 0 & 0\\ \hline
 \end{array}
 -
 \begin{array}{|r|r|r|}
  \multicolumn{3}{c}{} \\ \hline
  0 & 0 & 1\\ \hline
  1 & 0 & 0\\ \hline
  0 & 1 & 0\\ \hline
 \end{array}
 \end{equation}
 is a degree three move of the model (\ref{eq:cdem}). 
 Actually we easily see that row sums,  column sums and diagonal
 sums of $\Bz$ are all zeros. 
 By following the notation in (\ref{eq:moves}), $\bm{z}$ is written as
 \begin{equation}
  \label{eq:move-cell}
 \Bz = \left[
 \left\{
 (1,2), (2,3), (3,1) 
 \right\}
 ||  \left\{
 (3,2), (1,3), (2,1)
 \right\}
 \right].
 \end{equation}
 For this model only one move $\Bz$ forms a Markov basis (\cite{hara-takemura-yoshida-diagonal}).
\end{example}

Moves are used for steps of Markov chain Monte Carlo simulation within
each fiber.  
If we add or subtract a move $\Bz$ to $\Bx \in {\cal F}_{\Bb}$,
then $\Bx \pm \Bz \in {\cal F}_{\Bb}$ and we can move from $\Bx$ to
another state $\Bx+\Bz$ (or $\Bx - \Bz$) in the same fiber 
${\cal F}_\Bb$, as long as there is no negative element in 
$\Bx + \Bz$  (or $\Bx - \Bz$).

A finite set ${\cal M}$ of moves is called a {\em Markov basis} if for
every fiber the states become mutually accessible by the moves from
${\cal M}$. 
If we have a Markov basis, we can generate a Markov chain of contingency 
tables from any fiber whose stationary distribution is the conditional
hypergeometric distribution  (\citet{diaconis-sturmfels}).  
In this way Markov basis methodology enables us to evaluate a test
statistics based on the exact distribution.

\citet{dobra-2003bernoulli} showed that the decomposable model has a
Markov basis consisting of only degree two moves. 
Markov bases for some other log-affine model have been discussed in 
\citet{hty-jpaa}, \citet{hara-takemura-yoshida-diagonal} and 
\citet{th-thmc} etc.
In general, however it is not easy to obtain an exact list of Markov
basis for the log-affine model, even for the hierarchical model.
In hierarchical model \citet{dobra-sullivant} developed an algorithm to
compute a Markov basis recursively from Markov bases of the maximal
prime submodels corresponding to maximal compact components. 
In the next section we generalize the result to the HSM.

\subsection{Local computation of Markov basis of HSM}
Most of the arguments and the notations in this section follow 
those in \cite{dobra-sullivant}. 
For a subset $D \subset [m]$, denote $L(D) := L \cap F_D$.
Let $(A_1,A_2,S)$ be a decomposition of $L$ and define 
$V_1 := A_1 \cup S$ and $V_2 := A_2 \cup S$.  
Since $L$ is conformal to $\{F_{V_1}, F_{V_2}\}$, 
we note that ${\cal M}(L(V_1))$ and ${\cal M}(L(V_2))$ 
are marginal models
corresponding to $V_1$ and $V_2$, respectively.
Denote by $A_{V_1} = \{\Ba_{V_1}(\bm{i}_{V_1})\}_{\bm{i}_{V_1} \in {\cal I}_{V_1}}$ 
and $A_{V_2} = \{\Ba_{V_2}(\bm{i}_{V_2})\}_{\bm{i}_{V_2} \in {\cal I}_{V_2}}$ the 
configurations for the marginal models 
${\cal M}(L(V_1))$ and ${\cal M}(L(V_2))$,
where 
$\Ba_{V_1}(\bm{i}_{V_1})$ and $\Ba_{V_2}(\bm{i}_{V_2})$ denote column vectors of  
$A_{V_1}$ and $A_{V_2}$, respectively.
Noting that $\bm{i}_{V_1} = (\bm{i}_{A_1} \bm{i}_S)$ and $\bm{i}_{V_2} = (\bm{i}_S \bm{i}_{A_2})$, 
the configuration $A$ for ${\cal M}(L)$ is written by 
\[
 A = A_{V_1} \oplus_S A_{V_2} = 
 \{\Ba_{V_1}(\bm{i}_{A_1}\bm{i}_{S}) \oplus \Ba_{V_2}(\bm{i}_S \bm{i}_{A_2})\}_{\bm{i}_{A_1} \in
 {\cal I}_{A_1}, \bm{i}_S \in {\cal I}_S, \bm{i}_{A_2} \in {\cal I}_{A_2}}, 
\]
where 
\[
 \Ba_{V_1}(\bm{i}_{A_1}\bm{i}_{S}) \oplus \Ba_{V_2}(\bm{i}_S \bm{i}_{A_2}) = 
\left(
\begin{array}{c}
 \Ba_{V_1}(\bm{i}_{A_1}\bm{i}_{S})\\
 \Ba_{V_2}(\bm{i}_S \bm{i}_{A_2})
\end{array}
\right).
\]

Assume that ${\cal B}(V_1)$ and ${\cal B}(V_2)$ are  Markov bases for
${\cal M}(L(V_1))$ and ${\cal M}(L(V_2))$, respectively.  
Let 
$\bm{z}_1 = \{z_1(\bm{i}_{V_1})\}_{\bm{i}_{V_1} \in {\cal I}_{V_1}} \in
{\cal B}(V_1)$ and  
$\bm{z}_2 = \{z_2(\bm{i}_{V_2})\}_{\bm{i}_{V_2} \in {\cal I}_{V_2}} \in
{\cal B}(V_2)$. 
Since $S$ is saturated, the sufficient statistic $\Bb$ fixes $x(\Bi_S)$.
Hence we have 
\[
 \sum_{\bm{i}_{V_1 \setminus S} \in {\cal I}_{V_1 \setminus S}} z_1(\bm{i}_{V_1})
 = 0,  \quad  
 \sum_{\bm{i}_{V_2 \setminus S} \in {\cal I}_{V_2 \setminus S}} z_2(\bm{i}_{V_2})
 = 0.
\]
Then $\bm{z}_1$ and $\bm{z}_2$ can be written as
\begin{equation}
 \label{eq:z1}
  \bm{z}_1 = 
  [\{(\bm{i}^1_{A_1},\bm{i}^1_S),\ldots,(\bm{i}^d_{A_1},\bm{i}^d_S)\}||
  \{(\bm{j}^1_{A_1},\bm{j}^1_S),\ldots,(\bm{j}^d_{A_1},\bm{j}^d_S)\}], 
\end{equation}
\[
\bm{z}_2 = 
  [\{(\bm{i}^1_{S},\bm{i}^1_{A_2}),\ldots,(\bm{i}^d_S,\bm{i}^d_{A_2})\}||
  \{(\bm{j}^1_S,\bm{j}^1_{A_2}),\ldots,(\bm{j}^d_S,\bm{j}^d_{A_2})\}],
\]
respectively, where 
$\bm{i}^k_{A_1},\bm{j}^k_{A_1} \in {\cal I}_{A_1}$, 
$\bm{i}^k_S \in {\cal I}_S$ and 
$\bm{i}^k_{A_2},\bm{j}^k_{A_2} \in {\cal I}_{A_2}$
for $k=1,\ldots,d$. 

\begin{definition}[\citet{dobra-sullivant}]
 \label{def:dobra-sullivant}
 Define $\Bz_1 \in {\cal B}(V_1)$ as in (\ref{eq:z1}).
 Let 
 $\bm{\eta} := \{\bm{i}^1_{A_2},\ldots,\bm{i}^d_{A_2}\} \in 
 {\cal I}_{A_2} \times \cdots \times {\cal I}_{A_2}$.
 Define $\bm{z}_1^{\bm{k}}$ by 
 \[
 \bm{z}_1^{\bm{\eta}} := 
 [\{
 (\bm{i}^1_{A_1},\bm{i}^1_S,\bm{i}^1_{A_2}),
 \ldots,
 (\bm{i}^d_{A_1},\bm{i}^d_S,\bm{i}^d_{A_2})
 \}||
 \{
 (\bm{j}^1_{A_1},\bm{j}^1_S,\bm{i}^1_{A_2}),\ldots,
 (\bm{j}^d_{A_1},\bm{j}^d_S,\bm{i}^d_{A_2})
 \}].  
 \]
 Then we define 
 $\mathrm{Ext}({\cal B}(V_1) \rightarrow L)$ by 
 \[
 \mathrm{Ext}({\cal B}(V_1) \rightarrow L) := 
 \{\bm{z}_1^{\bm{\eta}} \mid \bm{\eta} \in {\cal I}_{A_2} \times \cdots
 \times {\cal I}_{A_2} \}. 
 \]
\end{definition}

In the same way as Lemma 5.4 in \citet{dobra-sullivant} 
we can obtain the following lemma.

\begin{lemma}
 \label{lemma:move-Ext}
 Suppose that $z_1 \in {\cal B}(V_1)$ as in (\ref{eq:z1}). 
 Then $\mathrm{Ext}({\cal B}(V_1) \rightarrow L)$  is the set of moves for
 $L$.  
\end{lemma}

\begin{proof}
 Let $\Bz \in \mathrm{Ext}({\cal B}(V_1) \rightarrow L)$. 
 Then we have 
 \[
  A \Bz = 
 \left(
 \begin{array}{c}
  \sum_{\bm{i}_{V_1} \in {\cal I}_{V_1}} \Ba_{V_1}(\bm{i}_{V_1})
   z_{V_1}(\bm{i}_{V_1})\\ 
  \sum_{\bm{i}_{V_2} \in {\cal I}_{V_2}} \Ba_{V_2}(\bm{i}_{V_2})
   z_{V_2}(\bm{i}_{V_2})\\ 
 \end{array}
 \right), 
 \]
 where 
 \[
 z_{V_1}(\bm{i}_{V_1}) = 
 \sum_{\bm{i}_{V_1^C} \in {\cal I}_{V_1^C}} z(\bm{i}), \quad 
 z_{V_2}(\bm{i}_{V_2}) = 
 \sum_{\bm{i}_{V_2^C} \in {\cal I}_{V_2^C}} z(\bm{i}).
 \]
 Since $z_{V_1}(\bm{i}_{V_1}) = z_1(\bm{i}_{V_1})$ and 
 $z_1 \in {\cal B}(V_1)$,  
 $\sum_{\bm{i}_{V_1} \in {\cal I}_{V_1}} \Ba_{V_1}(\bm{i}_{V_1})
 z_{V_1}(\bm{i}_{V_1})=0$. 
 From Definition \ref{def:dobra-sullivant}, 
 $z_{V_2}(\bm{i}_{V_2})=0$ for all $\bm{i}_{V_2} \in {\cal I}_{V_2}$.
 Hence $A\Bz = 0$.
\end{proof}

\begin{example}
 Consider a $3 \times 3 \times 3$ model in the class (\ref{eq:3way-ci-sub}), 
 \begin{equation}
  \label{eq:3way-cdem}
 \log p(\Bi) =\alpha(i_1) + \beta(i_2) + \gamma(i_3) 
  + \delta \phi(\Bi_{\{1,2\}}) + \delta' \phi(\Bi_{\{2,3\}}), 
 \end{equation}
 where $\phi(\cdot)$ is defined as in (\ref{eq:phi-cdem}). 
 The sufficient statistic for this model is the set of 
 one dimensional marginals $x(i_k)$, $i_k \in {\cal I}_k$, $k=1,2,3$ and
 two dimensional  diagonal sums  $\sum_{\Bi : i_1 = i_2} x(\Bi)$, 
 $\sum_{\Bi : i_2 = i_3} x(\Bi)$. 

 As discussed in Example \ref{ex:HSM},  this model is an HSM of
 (\ref{eq:cond-indep}). 
 Hence we can set $V_1 = \{1,2\}$ and $V_2 = \{2,3\}$ and 
 $L(V_i)=L\cap F_{V_i}$, $i=1,2$, are both 
 $3 \times 3$ common diagonal effect models (\ref{eq:cdem}). 

 Let $\Bz_1 := \Bz$ in (\ref{eq:move-cell}). 
 As mentioned in Example \ref{cdem}, $\Bz_1$ forms a Markov basis for
 the model (\ref{eq:cdem}), that is,  
 ${\cal B}(V_1) = \{\Bz_1\}$. 
 We  see that $\Bz_1$ is written in the form (\ref{eq:z1}). 
 Let $\bm{\eta} := (i_3,i_3',i_3'')$.
 Then $\Bz_1^{\bm{\eta}}$ is written by
 \[
 \Bz_1^{\bm{\eta}} = \left[
 \left\{
 (1,2,i_3), (2,3,i_3'), (3,1,i_3'') 
 \right\}
 ||  \left\{
 (3,2,i_3), (1,3,i_3'), (2,1,i_3'')
 \right\}
\right].
 \]
 When $\bm{\eta} = (1,2,3)$, $\Bz_1^{\eta}$ is written in array expression as in
 (\ref{eq:move-array}) by 
 \[
 \Bz_1^{\bm{\eta}} = 
 \begin{array}{ccc}
  \begin{array}{|r|r|r|}
    \multicolumn{3}{c}{i_3 = 1} \\ \hline
    0 & 1 & 0\\ \hline
    0 & 0 & 0\\ \hline
    0 & -1 & 0\\ \hline
  \end{array}, \; & 
  \begin{array}{|r|r|r|}
    \multicolumn{3}{c}{i_3 = 2} \\ \hline
    0 & 0 & -1\\ \hline
    0 & 0 & 1\\ \hline
    0 & 0 & 0\\ \hline
  \end{array}, \; & 
  \begin{array}{|r|r|r|}
    \multicolumn{3}{c}{i_3 = 3} \\ \hline
    0 & 0 & 0\\ \hline
    -1 & 0 & 0\\ \hline
    1 & 0 & 0\\ \hline
  \end{array}
 \end{array}.
 \]
 We easily see that 
 one dimensional marginals and two dimensional  diagonal sums of 
 $\Bz_1^{\bm{\eta}}$ are all zeros and hence that 
 $\Bz_1^{\bm{\eta}}$ is a move for (\ref{eq:3way-cdem}). 
 $\mathrm{Ext}({\cal B}(V_1) \rightarrow L)$ is 
 \[
 \mathrm{Ext}({\cal B}(V_1) \rightarrow L) = 
 \{
  \Bz_1^{\bm{\eta}} \mid i_3,i_3',i_3''\in \{1,2,3\}
 \}.
 \]
\end{example}

Consider a decomposable model ${\cal M}(H_\Delta)$ such that $\red\Delta =
\{V_1,V_2\}$.  
\citet{dobra-2003bernoulli} showed that the set of all degree two moves
\[
\Bz_{V_1,V_2} = \left[
\left\{(\Bi^1_{A_1}, \Bi^1_{S}, \Bi^1_{A_2}), (\Bi^2_{A_1}, \Bi^2_{S}, \Bi^2_{A_2}) \right\} || 
\left\{(\Bi^1_{A_1}, \Bi^1_{S}, \Bi^2_{A_2}), (\Bi^2_{A_1}, \Bi^2_{S}, \Bi^1_{A_2}) \right\} 
\right], 
\]
where 
$\bm{i}^k_{A_1} \in {\cal I}_{A_1}$, $\bm{i}^d_S \in {\cal I}_S$ and
$\bm{i}^k_{A_2} \in {\cal I}_{A_2}$ for $k=1,2$,  
forms a Markov basis and denote it by ${\cal B}_{V_1,V_2}$.

\begin{theorem}
 \label{thm:dobra-sullivant}
 Let ${\cal B}(V_1)$ and ${\cal B}(V_2)$ be Markov bases for ${\cal M}(L(V_1))$ and
 ${\cal M}(L(V_2))$, respectively. 
 Then 
 \begin{equation}
  \label{eq:dobra-sullivant-recursive}
 {\cal B} := \mathrm{Ext}({\cal B}(V_1) \rightarrow L) \cup 
 \mathrm{Ext}({\cal B}(V_2) \rightarrow L) \cup 
 {\cal B}_{V_1,V_2}
 \end{equation}
 is a Markov basis for ${\cal M}(L)$. 
\end{theorem}

We can prove the theorem in the same way as Theorem 5.6 in
\citet{dobra-sullivant}. 
Suppose that ${\cal M}(L)$ is 
an HSM of $H_{\cal H}$.
Then Theorem \ref{thm:dobra-sullivant} implies that a Markov basis for
$L$ is obtained from ${\cal B}(C)$, $C \in {\cal H}$, by recursively
using (\ref{eq:dobra-sullivant-recursive}). 
This shows that the computation of a Markov basis can be localized 
according to reducible submodels corresponding to maximal extended
compact components of $L$. 

Concerning Markov bases of the split model of Section \ref{sec:split} we state the following lemma.

\begin{lemma}
 \label{lemma:split-Markov-1}
 With the same notation as in Lemma \ref{lemma:split-Markov},  a
 Markov basis of the model associated with the subspace $\sum_{\lambda}N_{\lambda}$  
 is given by union of Markov bases of ${\cal M}(N_{\lambda})$. 
\end{lemma}

\section{Examples}
\label{sec:examples}

In this section we give several applications of conditional tests of 
HSMs by using Markov bases. 
In Section \ref{subsec:pewee} 
we discuss conditional tests for models of multiway tables with
structural zeros. 
In Section \ref{subsec:wam} we present an example of a split model.
The models in this section are relatively small and intended to
illustrate the notions of this paper, rather than being  examples of large scale
data analyses.

\subsection{Conditional tests for models with structural zeros}
\label{subsec:pewee}

\begin{table}[ht]
 \caption{Triples of phrases in a song sequence
 of a wood pewee, with repeats deleted.}
 \label{tab:wood-pewee}
 \begin{center}
  {\small
  \begin{tabular}{cccccc}
   \hline
   && \multicolumn{4}{c}{Third place}\\
   First place& Second place& A& B& C& D\\
   \cline{3-6}
   A& A& ---& ---& ---& ---\\
   & B& 19& ---& 2& 2\\
   & C& 2& 26& ---& 0\\
   & D& 12& 5& 0& ---\\
   B& A& ---& 9& 6& 12\\
   & B& ---& ---& ---& ---\\
   & C& 24& 1& ---& 1\\
   & D& 1& 2& 0& ---\\
   C& A& ---& 4& 22& 0\\
   & B& 3& ---& 22& 0\\
   & C& ---& ---& ---& ---\\
   & D& 1& 0& 0& ---\\
   D& A& ---& 11& 0& 4\\
   & B& 5& ---& 1& 1\\
   & C& 0& 0& ---& 0\\
   & D& ---& ---& ---& ---\\
   \hline
   \multicolumn{6}{l}{{\it Source:} \citet{Craig}}
  \end{tabular}
  }
 \end{center}
\end{table}

Table \ref{tab:wood-pewee} is the data on song sequence of a wood pewee
in Section 7.5.2 of \citet{bishop-fienberg-holland}. 
The wood pewee has a repertoire of four distinctive phrases.
The observed data consists of 198 triplets of consecutive phrases
$(i,j,k)\in \{1,2,3,4\}^3$. 
It is a $4\times 4\times 4$ contingency table with the cells of the form
$(i,i,k)$ and $(i,j,j)$ being 
structural zeros. 
As discussed in \cite{chatfield-lemon}, we consider this sequence as a
Markov chain. 
The main interest is the order of the chain. 
As an example of conditional tests for the model with structural zeros, 
we consider the goodness-of-fit test of two Markov chain models of first
order for this data. 
\citet{aoki-takemura-2005jscs} provided a complete description of Markov
basis for the quasi-independence model for two-way tables and proposed
conditional test by using the Markov basis.
However its extension to the model for multiway tables has not yet been
studied. 

First we consider the model discussed by \citet{bishop-fienberg-holland} 
for this data, 
\begin{equation}
\label{eq:pewee}
p_{ijk}= 1_{\{i\neq j\}} e^{a_{ij}}  1_{\{j\neq k\}}  e^{b_{jk}}, 
\end{equation}
where $a_{ij}$ and $b_{jk}$ are free parameters.
With some abuse of notation (\ref{eq:pewee}) can be written as
\begin{align}
\label{eq:pewee2}
 \log p_{ijk}
 \ =\  a_{ij}1_{\{i\neq j\}}+(-\infty)1_{\{i=j\}}
 +b_{jk}1_{\{j\neq k\}}+(-\infty)1_{\{j=k\}}.
\end{align}
We note that this model is also in the class (\ref{eq:3way-ci-sub}). 
The probability function $\{p_{ijk}\}$ satisfies the condition
$p_{iik}=0$ and $p_{ijj}=0$, or equivalently,
$\log p_{iik}=-\infty$ and $\log p_{ijj}=-\infty$.
Hence $\{\log p_{ijk}\}$ is not an element of $V=\R^{4\times 4\times 4}$.
However we can replace $V$ by $R^{|\bar {\cal I}|}$, where
\[
\bar {\cal I}= {\cal I} \setminus \big(\{ (i,i,j), i,j\in [4] \}  \cup 
\{ (i,j,j), i,j\in [4] \}  \big), 
\]
and consider log-affine models of $R^{|\bar {\cal I}|}$. Formally it is more convenient
to proceed with $V=\R^{4\times 4\times 4}$ allowing 
$\log p_{iik}=\log p_{ijj}=-\infty$.

\begin{figure}[t]
\begin{center}
\begin{tabular}{cc}
\includegraphics[width=0.4\textwidth]{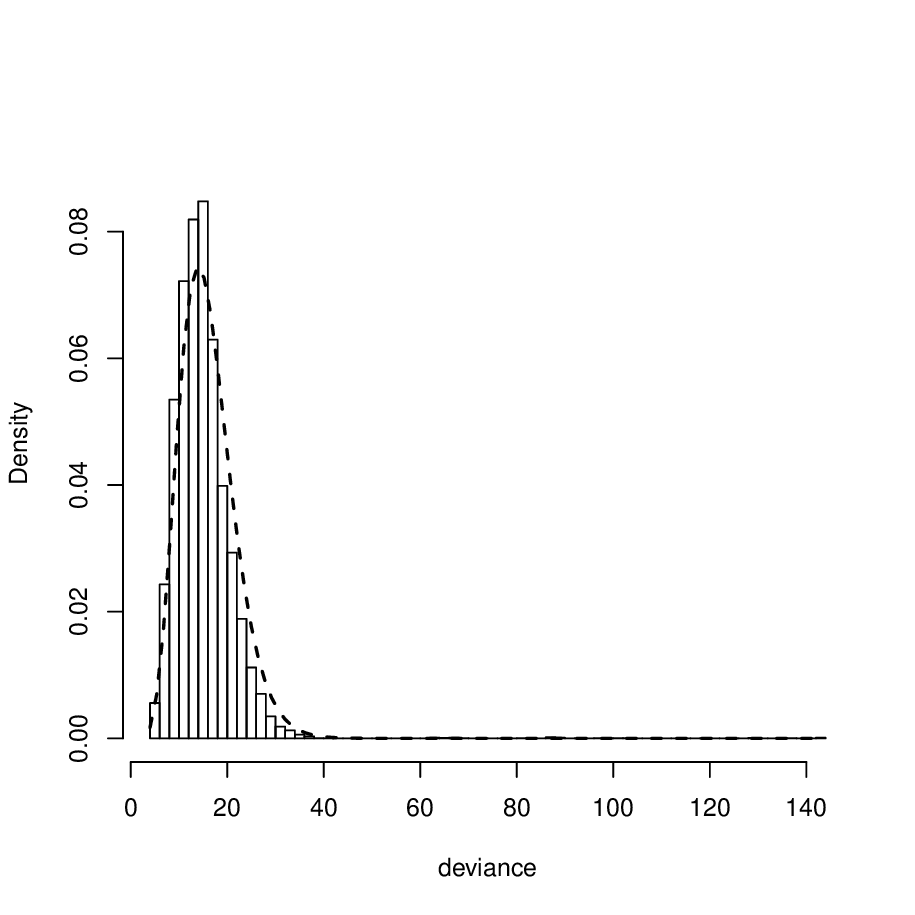}
 &
\includegraphics[width=0.4\textwidth]{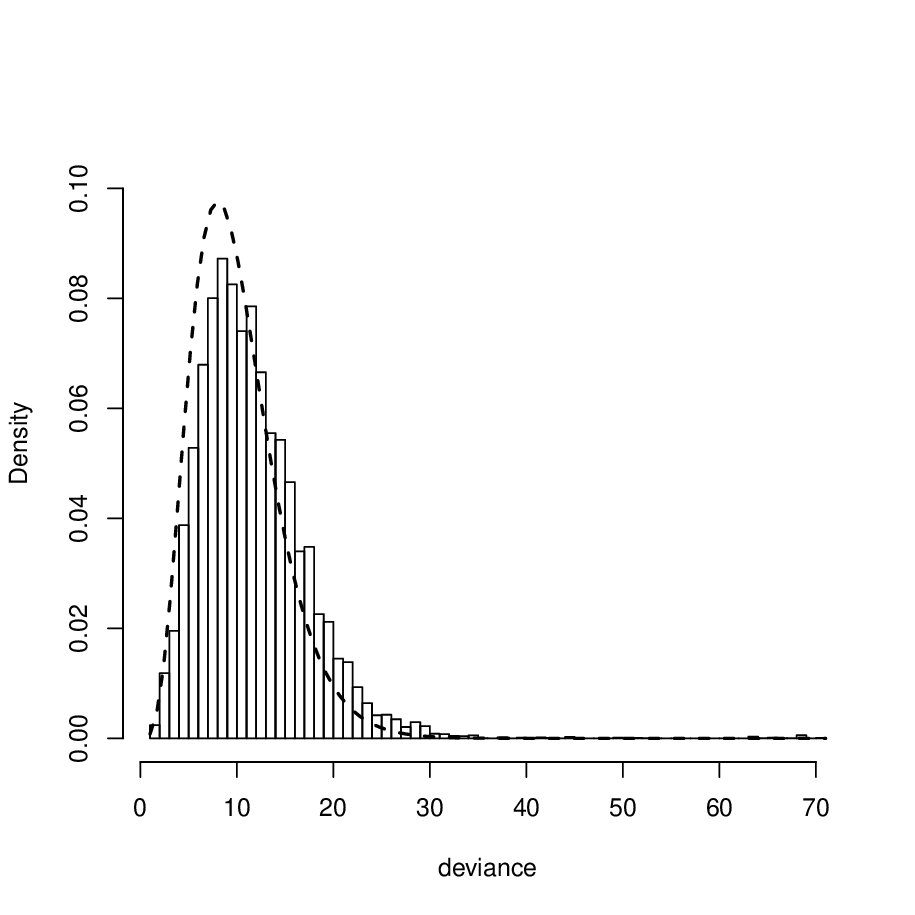}
 \\
 (a) Deviance of Model 1 ($G^2=142.4$).
 &
 (b) Deviance of Model 2 from Model 1 ($G^2=66.9$).
\end{tabular}
\end{center}
\caption{The empirical distribution and asymptotic distribution of 
 deviance $G^2$ for the wood pewee data.
 The degree of freedom
 is 16 and 10, respectively.
 The number of steps in the MCMC procedure is $10^5$.}
\label{fig:wood-pewee}
\end{figure}

We first consider the conditional independence model ${\cal M} (F_{\rm Model1})$,
where
\begin{align*}
 F_{\rm Model1}\ =\ F_{\{1,2\}}+F_{\{2,3\}},
\end{align*}
which corresponds to (\ref{eq:pewee}).
The MLE of this model is explicitly given by
\begin{align*}
 \hat{p}_{ijk}
 \ =\ \frac{x_{ij+}x_{+jk}}{n x_{+j+}}
 \ =\ \frac{x_{ij+}1_{\{i\neq j\}}x_{+jk}1_{\{j\neq k\}}}{n x_{+j+}}.
\end{align*}
A Markov basis of the model is
$\cB_{\rm Model1}=\cB_{\{1,2\},\{2,3\}}$
(see Theorem~\ref{thm:dobra-sullivant} for the notation).
An experimental result that compares the saturated model
and Model 1 is given in Figure~\ref{fig:wood-pewee}.
Both the asymptotic and experimental estimates of the p-value
are almost zero.

Although Model~1 does not fit the data,
we proceed to consider a submodel of Model~1 for theoretical interest.
Let
\begin{align*}
 F_{\rm model2}\ =\ 
 \left\{
 \alpha_i+\beta_j+\gamma_k+\phi_i1_{\{i=j\}}+\psi_j1_{\{j=k\}}
 \right\}.
\end{align*}
${\cal M}(F_{\rm model2})$ 
is 
an HSM of $F_{\{1,2\}}+F_{\{2,3\}}$.
It represents a quasi-independence model for the three-way table.
The MLE of the model is
\begin{align*}
 \hat{p}_{ijk}
 \ =\ \frac{\hat{p}^{(1)}_{ij}\hat{p}^{(2)}_{jk}}{x_{+j+}/n},
\end{align*}
where $\hat{p}^{(1)}_{ij}$ and $\hat{p}^{(2)}_{jk}$
are the MLE of the 2-way quasi-independence models
with the diagonal structural zeros, that is,
\begin{align*}
& \hat{p}^{(1)}_{ij}=e^{\hat{\alpha}_i}e^{\hat{\beta}_j}1_{\{i\neq j\}},
\quad \hat{p}^{(1)}_{i+}=x_{i++}/n,
\quad \hat{p}^{(1)}_{+j}=x_{+j+}/n,
 \\
& \hat{p}^{(2)}_{jk}=e^{\hat{\beta}'_j}e^{\hat{\gamma}_k}1_{\{j\neq k\}},
\quad \hat{p}^{(2)}_{j+}=x_{+j+}/n,
\quad \hat{p}^{(2)}_{+k}=x_{++k}/n,
\end{align*}
where $\hat{\beta}_j$ and $\hat{\beta}'_j$ are different in general
as discussed in Example~\ref{ex:3way}.
They are computed by the iterative proportional fitting method.
By Theorem~\ref{thm:dobra-sullivant}, a Markov basis is given by
\begin{align*}
 \cB_{\rm Model2}\ =\ \cB_{\{1,2\},\{2,3\}}
 \cup \mathrm{Ext}(\cB(\{1,2\})\to V)
 \cup \mathrm{Ext}(\cB(\{2,3\})\to V)
\end{align*}
where $\cB(\{1,2\})$ and $\cB(\{2,3\})$
are the Markov bases of the 2-way quasi-independence model with
structural zeros obtained by \citet{aoki-takemura-2005jscs}.
An experimental result that compares the Model 1
and Model 2 is given in Figure~\ref{fig:wood-pewee}.
These results show that we can conclude the chain is at least of second
order. 

In this way we can perform conditional test for the models of multiway
tables with conditional zeros.

\subsection{Conditional test for the split model}
\label{subsec:wam}

In this section we give an example of conditional test of the split
model. 
Here we deal with a real data called {\em women and mathematics (wam)
data} used in \citet{Hojsgaard}. 
The data is shown in Table~\ref{tab:wam}.
The data consists of the following six factors:
(1) Attendance in math lectures (attended=1, not=2),
(2) Sex (female=1, male=2),
(3) School type (suburban=1, urban=2),
(4) Agree in statement ``I'll need mathematics in my future work'' (agree=1, disagree=2),
(5) Subject preference (math-science=1, liberal arts=2)
and (6) Future plans (college=1, job=2).
We consider two models \citet{Hojsgaard} treated.
The first model is a decomposable model ${\cal M}(F_{\rm model1})$ 
\begin{align*}
 F_{\rm Model1}\ =\ F_{\{1,2,3,5\}} + F_{\{2,3,4,5\}} + F_{\{3,4,5,6\}}.
\end{align*}
By Theorem~\ref{thm:dobra-sullivant},
a Markov basis of this model is given by
\begin{align*}
 \cB_{\rm Model1}\ =\ \cB_{\{1,2,3,5\},\{2,3,4,5,6\}}\cup\cB_{\{1,2,3,4,5\},\{3,4,5,6\}}.
\end{align*}

\begin{table}[ht]
 \caption{Survey data concerning the attitudes of
 high-school students in New Jersey towards mathematics.}
 \label{tab:wam}
 \begin{center}
  {\small
  \begin{tabular}{llcccccccc}
   \hline
   School& & \multicolumn{4}{l}{Suburban school}& \multicolumn{4}{l}{Urban school}\\
   \hline
   Sex& & \multicolumn{2}{l}{Female}& \multicolumn{2}{l}{Male}& \multicolumn{2}{l}{Female}& \multicolumn{2}{l}{Male}\\
   \cline{2-4} \cline{5-6} \cline{7-8} \cline{9-10}
   Plans& Preference& Attend& Not& Attend& Not& Attend& Not& Attend& Not\\
   \hline
   College& Math-sciences&&&&&&\\
   & Agree&     37&  27&  51&  48&  51&  55&  109&  86\\
   & Disagree&  16&  11&  10&  19&  24&  28&  21&  25\\
   & Liberal arts&&&&&&\\
   & Agree&     16&  15&   7&   6&  32&  34&  30&  31\\
   & Disagree&  12&  24&  13&   7&  55&  39&  26&  19\\
   Job& Math-sciences&&&&&&\\
   & Agree&     10&   8&  12&  15&   2&   1&   9&   5\\
   & Disagree&   9&   4&   8&   9&   8&   9&   4&   5\\
   & Liberal arts&&&&&&\\
   & Agree&      7&  10&   7&   3&   5&   2&   1&   3\\
   & Disagree&   8&   4&   6&   4&  10&   9&   3&   6\\
   \hline
   \multicolumn{10}{l}{{\it Source:} \citet{Fowlkes}}\\
  \end{tabular}
  }
 \end{center}
\end{table}

The second model is a split model ${\cal M}(F_{\rm model2})$ 
\begin{align*}
 F_{\rm Model2}\ =\ 
 F_{\{1,2,3,5\}} + F_{\{2,5\}}^{j_3=1} + F_{\{4,5\}}^{j_3=1}
 + F_{\{2,4,5\}}^{j_3=2} + F_{\{3,4,5,6\}}.
\end{align*}
This model is indeed a split model (of degree one) with
\begin{align*}
 & \cC\ =\ \{\{1,2,3,5\},\{2,3,4,5\},\{3,4,5,6\}\},
 \\
 & Z(\{1,2,3,5\})\ =\ \emptyset,
 \quad \cC_{\{1,2,3,5\}}^{\Bj_{\emptyset}}\ =\ \{\{1,2,3,5\}\},
 \\
 & Z(\{2,3,4,5\})\ =\ \{3\},
 \quad \cC_{\{2,3,4,5\}}^{j_3=1}\ =\ \{\{2,5\},\{4,5\}\},
 \quad \cC_{\{2,3,4,5\}}^{j_3=2}\ =\ \{\{2,4,5\}\},
 \\
 & Z(\{3,4,5,6\})\ =\ \emptyset,
 \quad \cC_{\{3,4,5,6\}}^{\Bj_{\emptyset}}\ =\ \{\{3,4,5,6\}\}.
\end{align*}
The condition (\ref{eqn:split-saturated}) is easily checked.
The MLE is calculated if one decomposes the table
into those for $j_3=1$ and $j_3=2$ and then calculates
the MLE separately (Lemma~\ref{lemma:split-Markov}).
By Theorem~\ref{thm:dobra-sullivant} and Lemma~\ref{lemma:split-Markov-1},
a Markov basis of this model is
\begin{align*}
 \cB_{\rm Model2}\ =\ 
 \cB_{\{1,2,5\},\{4,5,6\}}^{j_3=1}
 \cup \cB_{\{1,2,3,5\},\{2,3,4,5,6\}}
 \cup \cB_{\{1,2,3,4,5\},\{3,4,5,6\}},
\end{align*}
where we put
$\cB_{\{1,2,5\},\{4,5,6\}}^{j_3=1}
=\cB_{\{1,2,5\},\{4,5,6\}}\cap F^{j_3=1}$.

\begin{figure}[ht]
 \begin{center}
  \includegraphics[width=0.4\textwidth]{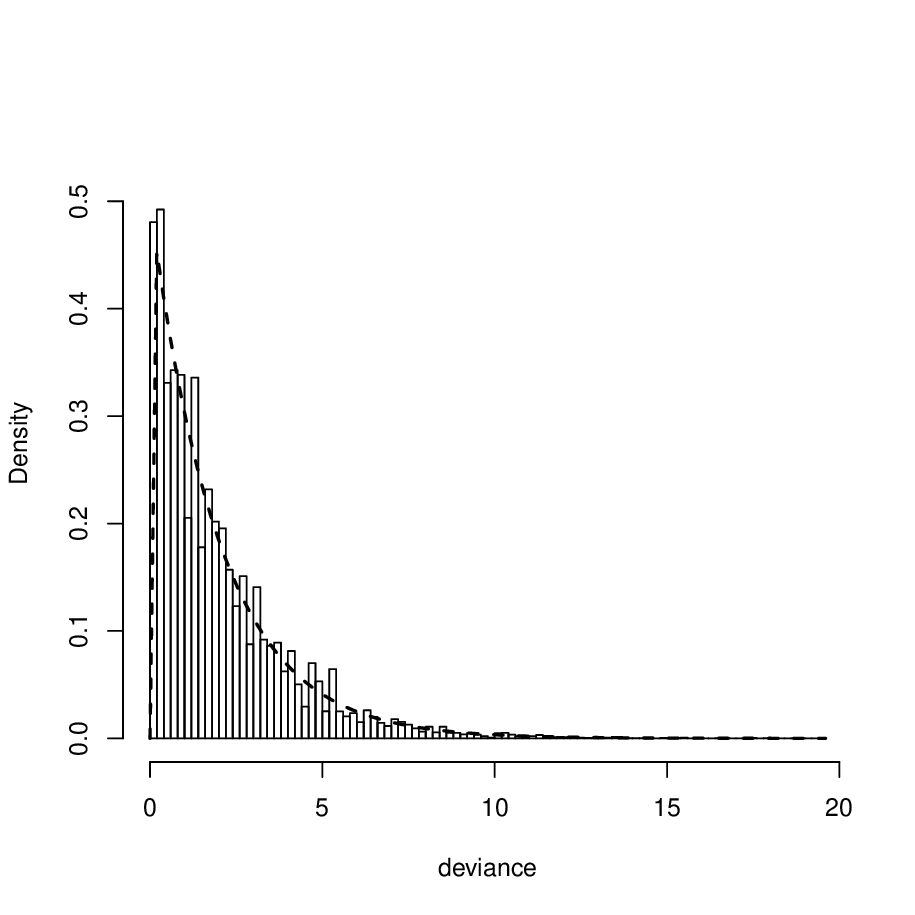}
 \end{center}
 \caption{The empirical and asymptotic distributions
 of the deviance of Model 2 from Model 1.}
 \label{fig:wam}
\end{figure}

We calculate the p-value of the deviance
of Model 2 from Model 1 by the MCMC method.
The number of steps in the MCMC procedure is $10^5$.
The result is as follows.
\begin{center}
\begin{tabular}{cccc}
 \hline
 Deviance & df & p-value (asymptotic) & p-value (MCMC)\\
 \hline
 1.851 & 2 & 0.396 & 0.399$\pm$0.012\\
 \hline
\end{tabular}
\end{center}
The confidence interval of the p-value
is computed on the basis of the batch-means method.
The empirical distribution and asymptotic distribution
of the deviance are given in Figure~\ref{fig:wam}.
In this way we can perform conditional test for the split model.

\section{Concluding remarks}
\label{sec:remarks}
We proposed a hierarchical subspace model, by defining
the notion of conformality of linear subspaces to a given hierarchical
model.  The notion of 
an HSM gives a modeling
strategy of multiway tables and unifies various models of
interaction effects in the literature.  
We illustrated 
our modeling strategy with some data sets.
As a referee pointed out, our approach is novel in the sense that
the localization properties are described not only
by means of graph-theoretical criteria but also using the properties
of the linear subspaces encoding these models.

In this paper we only considered log-affine model.  Note that there
are some nonlinear models of interaction terms for two-way tables, such
as the RC association model. It seems clear that we can separately fit
a nonlinear model to each maximal compact component of a hierarchical model, as long
as the models for dividers are saturated. However conformality of
a general nonlinear model with respect to a given hierarchical model has to be 
carefully defined and this is left to our future study.

The separation by dividers are closely related to the notion of collapsibility
(e.g.\ \citet{asmussen1983}) of hierarchical models.  Localization of statistical
inference to the marginal table of a maximal extended compact component seems to correspond
to the collapsibility to the component.  Also Theorem 1 suggests the effectiveness 
of using mixed parameterization for contingency tables, i.e., we fit log-linear models for 
maximal extended compact components and connect them 
by marginal probabilities as in (\ref{eq:mle}).
Furthermore  our results for Markov bases for HSMs
are 
closely related to those of \citet{sullivant-toric-fiber}.  \citet{sullivant-toric-fiber} 
is more concerned with Markov bases for models with 
latent variables and marginalization of latent variables.
Collapsibility and marginalization properties of 
HSM require further investigation.

In the computation of the MLE for the hierarchical models, 
it is known that the algorithm can be localized into the marginal
tables of maximal cliques for chordal extension of the simplicial
complex associated with the model, which is smaller than maximal compact
component (e.g.\ \citet{Badsberg-Malvestuto}). 
By using the notion of ambient hierarchical model discussed in Section 
\ref{subsec:closure}, it may be possible to localize the inference to
smaller units than maximal extended compact component also in the 
HSMs. 

Another important question on hierarchical subspace model is the
necessity of saturation of the model for dividers.  Saturation of the model for dividers
is a sufficient condition for localization of statistical inference, but it may not be
a necessary condition.  There may exist some important models, for
which statistical inferences can be localized to extended compact
components without the requirement of saturation of dividers.  This question also needs
a careful investigation.

\ \\
\noindent {\bf Acknowledgments.}\ 
The authors are grateful to three anonymous referees for constructive  
and detailed comments. 



\nocite{th-thmc}
\bibliographystyle{plainnat}
\bibliography{HST-HSM}

\end{document}